\documentclass[12pt]{amsart}

  \usepackage{amssymb}
  \usepackage{amsmath}
  \usepackage{latexsym}
  \usepackage{amsthm}
  \usepackage[all]{xy}
  \usepackage{amsbsy}

  \makeatletter

  \def\thmhead@plain#1#2#3{%
  \thmname{#1}\thmnumber{\@ifnotempty{#1}{ }#2}%
  \thmnote{ \the\thm@notefont(#3)}}

  \let\thmhead\thmhead@plain

  \def\swappedhead#1#2#3{%
  \thmnumber{#2}\thmname{\@ifnotempty{#2}{. }#1}%
  \thmnote{ \the\thm@notefont(#3)}}

  \makeatother

 \theoremstyle{definition} 

 \newtheorem{definition}{Definition}[section]
 \newtheorem{remark}[definition]{Remark}
 \newtheorem{example}[definition]{Example}
 
 \newtheorem{question}[definition]{Question}

 \theoremstyle{plain}      

 \newtheorem{proposition}[definition]{Proposition}
 \newtheorem{theorem}[definition]{Theorem}
 \newtheorem{corollary}[definition]{Corollary}
 \newtheorem{lemma}[definition]{Lemma}

\renewcommand{\colon}{ : }

\newcommand{\RRaro}{\Longrightarrow}

\newcommand{\raro}{\rightarrow}

\newcommand{\cF}{{\mathcal F}}
\newcommand{\cS}{{\mathcal S}}
\newcommand{\cR}{{\mathcal R}}
\newcommand{\cL}{{\mathcal L}}

\newcommand{\bH}{{\bar H}}
\newcommand{\bh}{{\bar h}}
\newcommand{\bchi}{{\bar \chi}}
\newcommand{\bT}{{\bar T}}
\newcommand{\bu}{{\bar u}}
\newcommand{\bI}{{\bar I}}
\newcommand{\bE}{{\bar E}}
\newcommand{\baf}{{\bar f}}
\newcommand{\bxi}{{\bar \xi}}
\newcommand{\bC}{{\bar C}}
\newcommand{\bpsi}{{\bar \psi}}
\newcommand{\bphi}{{\bar \phi}}
\newcommand{\bcS}{{\bar \cS}}
\newcommand{\bXi}{{\bar \Xi}}
\newcommand{\bgk}{{\bar \gk}}
\newcommand{\bgga}{{\bar \gga}}
\newcommand{\bgb}{{\bar \gb}}

\newcommand{\bldphi}{\boldsymbol{\phi}}
\newcommand{\bldrho}{\boldsymbol{\rho}}
\newcommand{\bldga}{\boldsymbol{\alpha}}
\newcommand{\bldgb}{\boldsymbol{\beta}}
\newcommand{\bldgga}{\boldsymbol{\gamma}}
\newcommand{\bldpi}{\boldsymbol{\pi}}
\newcommand{\bldpsi}{\boldsymbol{\psi}}
\newcommand{\bldd}{\mathbf{d}}
\newcommand{\bldT}{\mathbf{T}}
\newcommand{\bldl}{\boldsymbol{(}}
\newcommand{\bldr}{\boldsymbol{)}}
\newcommand{\bldlTr}{\bldl \bldT \bldr}
\newcommand{\bldG}{\mathbf{G}}
\newcommand{\bldS}{\mathbf{S}}
\newcommand{\bldP}{\mathbf{P}}
\newcommand{\bldGSP}{\bldG \bldS \bldP}
\newcommand{\bldu}{\mathbf{u}}
\newcommand{\bldp}{\mathbf{p}}

\newcommand{\hpsi}{{\hat \psi}}
\newcommand{\hphi}{{\hat \phi}}
\newcommand{\hH}{{\hat H}}

\newcommand{\ga}{\alpha}
\newcommand{\gga}{\gamma}
\newcommand{\gG}{\Gamma}
\newcommand{\gb}{\beta}
\newcommand{\gd}{\delta}
\newcommand{\gD}{\Delta}
\newcommand{\gk}{\kappa}

\newcommand{\gep}{\varepsilon}

\newcommand{\gf}{\varphi}

\newcommand{\gs}{\sigma}

\newcommand{\gw}{\omega}

\newcommand{\gz}{\zeta}

\newcommand{\N}{{\mathbb{N}}}

\newcommand{\Z}{{\mathbb{Z}}}
\newcommand{\R}{{\mathbb{R}}}

\newcommand{\pn}{^{(n)}}
\newcommand{\pinf}{^{(\infty)}}
\newcommand{\po}[1]{^{(#1)}}
\newcommand{\ti}[1]{\tilde{#1}}
\newcommand{\Hom}{{\mathcal Hom}}
\newcommand{\Ext}{{\mathcal Ext}^d(I,\rho)}
\newcommand{\Ad}{{\mathcal A}_d}
\newcommand{\UE}{\mathcal {UE}(\rho)}
\newcommand{\SUE}{\mathcal {SUE}(\rho)}

\begin{document}

\title{ On a class of one-sided Markov shifts}

\author{Ben-Zion Rubshtein}

\markboth{Ben Zion Rubshtein}{On one-sided Markov shifts }

\address{Address: {\it Ben Zion Rubshtein,
Dept. of Mathematics, Ben-Gurion University of the Negev,
Beer-Sheva, 84105, Israel.}}

\address{ E-mail: {\it benzion@math.bgu.ac.il}}


\begin{abstract}

We study one-sided Markov shifts, corresponding to positively
recurrent Markov chains with countable (finite or infinite) state
spaces.
The following classification problem is considered: when
two one-sided  Markov shifts are isomorphic up to a measure
preserving isomorphism ?
In this paper we solve the problem for
the class of $\rho$-uniform (or finitely $\rho$-Bernoulli)
one-sided Markov shifts considered in \cite{Ru$_6$}.

We show that every ergodic $\rho$-uniform Markov shift $T$ can be
represented in a canonical form $T = T_G $ by means of a canonical
(uniquely determined by $T$) stochastic graph $G$. In the
canonical form, two such shifts $T_{G_1}$ and $T_{G_2}$ are
isomorphic if and only if their canonical stochastic graphs $G_1$
and $G_2$ are isomorphic.

\end{abstract}

\maketitle


\section{Introduction}\label{Intro}

\bigskip

In this paper we consider the classification problem for one-sided
Markov shifts with respect to measure preserving isomorphism. Let
$G$ be a finite or countable stochastic graph, i.e. a directed
graph, whose edges $g \in G$ are equipped with positive weights
$p(g)$. The weights $p(g)$ determine transition probabilities of a
Markov chain on the discrete state space $G$. The corresponding
one-sided Markov shift $T_G$ acts on the space $(X_G,m_G)$, where
$ X_G = G^\N $ and $m_G$ is a stationary (probability) Markov
measure on $X_G$. We deal only with irreducible positively
recurrent Markov chains, so that such a Markov measure exists and
the shift $T_G$ is an ergodic endomorphism of the Lebesgue space
$(X_G,m_G)$.
The problem under consideration is : When for given
two stochastic graphs $G_1$ and $G_2$, does there exist an
isomorphism $\; \Phi \colon  X_{G_1} \raro  X_{G_2} \;$ such that
$\; m_{G_2} = m_{G_1} \circ \Phi^{-1} \;$ and
$\; \Phi \circ T_{G_1} \;=\; T_{G_2} \circ \Phi \;$.

It is obvious, that any (weight preserving) graph isomorphism
$\phi \colon G_1 \raro G_2$ generates such an isomorphism
$\Phi = \Phi_\phi$, but nonisomorphic graphs can generate the same
shift $T_G$.

\medskip

Recently J. Ashley, B. Marcus and S. Tuncel \cite{AsMaTu} solved the
classification problem for one-sided Markov shifts corresponding to
{\bf finite} Markov chains.
They used an approach which is based on the following important
fact: Two one-sided Markov shifts $T_{G_1}$ and $T_{G_2}$ (on
finite state spaces) are isomorphic iff there exists a common
extension $G$ of $G_1$ and $G_2$ by right resolving graph
homomorphisms of degree $1$. The result was proved implicitly in
\cite{BoTu}, where regular isomorphisms and right closing maps for
two-sided Markov shifts were studied (See also \cite{As},
\cite{KiMaTr}, \cite{Tr}, \cite{Ki} and references cited there)

\medskip

It should be noted that the classification problem for two-sided
shifts is quite different from the one-sided case.
Namely, any mixing two-sided Markov shift is isomorphic to the
Bernoulli shift with the same entropy \cite{FrOr} and two-sided
Bernoulli shifts are isomorphic iff they have the same entropy by
the Sinai-Ornstein theorem \cite{Si}-\cite{Or}.

On the other hand, let $T_\rho$ be the one-sided Bernoulli shift
with a discrete state space $(I,\rho)$, where $I$ is a finite or
countable set, $1 < |I| \leq \infty$, and $\rho = \{\rho_i\}_{i
\in I}$, $\sum_{i} \rho_i = 1$, $\rho_i > 0$.
The endomorphism $T_\rho$ acts as the one-sided shift on the
product space
$(X_\rho , m_\rho) = \prod_{n=1}^{\infty} (I, \rho)$.
Consider the measurable partition
$T_\rho^{-1} \gep = \{ T_\rho^{-1} x \;,\; x \in X_\rho \}$
generated by $T_\rho$ on $X_\rho$.
The partition admits an independent complement $\gd$, which is not
unique in general, but necessarily has the distribution $\rho$.
This implies that one-sided Bernoulli shifts $T_{\rho_1}$ and
$T_{\rho_2}$ are isomorphic iff the distributions $\rho_1$ and
$\rho_2$ coincide.

This simple observation motivates the following definition.
An endomorphism $T$ of a Lebesgue space $(X, m)$ is called
$\bldrho${\bf-uniform} (or {\bf finitely}
$\bldrho${\bf -Bernoulli} according to \cite{Ru$_6$}) if the
measurable partition
$T^{-1} \gep = \{ T^{-1} x \;,\;\; x \in X \}$ admits an
independent complement $\gd$ with $distr \; \gd = \rho$. We denote
by $\UE$ the class of all $\rho$-uniform endomorphisms.

Recall that the cofiltration $\xi (T)$ generated by an
endomorphism $T$ is the decreasing sequence
 $\{\xi_n\}_{n=1}^{\infty}$ of the measurable partitions
 $\xi_n = T^{-n} \gep$ of the space $X$ onto inverse images
 $T^{-n}x$.
 If two endomorphisms $T_1$ and $T_2$ are isomorphic, i.e.
there exists an isomorphism $\Phi$ such that
 $\; \Phi \circ T_1 = T_2 \circ \Phi \;$, then
 $\; \Phi (T_1^{-n}x) = T_2^{-n} (\Phi x) \;$ for almost all
 $x \in X$, i.e. $\Phi (T_1^{-n} \gep) = T_2^{-n} \gep$ for all
$n$.
 This means that the cofiltrations $\xi(T_1)$ and $\xi(T_2)$
are isomorphic.

  If $T \in \UE $, the cofiltration $\xi(T)$ is not necessarily
isomorphic to the {\bf standard} cofiltration $\xi(T_\rho)$,
generated by the Bernoulli shift $T_\rho$.
However, it is {\bf finitely isomorphic} to $\xi(T_\rho)$, i.e.
for every $n \in \N$ there exists an isomorphism $\Phi_n$ such
that $\Phi_n (T^{-k} \gep) =T_\rho^{-k} \gep$ for all $1 \leq k
\leq n$.

The isomorphism problem for $\rho$-uniform endomorphisms is
decomposed into the following two parts: When are the
cofiltrations $\xi(T_1)$ and $\xi(T_2)$ isomorphic? When are $T_1$
and $T_2$ isomorphic provided that $\xi(T_1) = \xi(T_2)$?

In particular, for given $ T \in \UE $: When is the cofiltration
$\xi(T)$ standard, i.e. isomorphic to $\xi(T_\rho)$ ? When are
$T_1$ and $T_2$ isomorphic provided that $\xi(T_1) = \xi(T_2)$ ?

All these problems are quite nontrivial even in the dyadic case
$\rho = (\frac{1}{2},\frac{1}{2})$. Various classes of decreasing
sequences of measurable partitions were considered by A.M. Vershik
\cite{Ve$_1$}-\cite{Ve$_4$}, V.G. Vinokurov \cite{Vi}, A.M. Stepin

\cite{St} and by author
\cite{Ru$_1$},\cite{Ru$_2$},\cite{Ru$_4$}-\cite{Ru$_6$}.
A new remarkable progress in the theory is due to J. Feldman, D.J.
Rudolph, D. Heicklen and Ch. Hoffman (See \cite{FeR}, \cite{HeHo},
\cite{HeHoR}, \cite{Ho}, \cite{HoR}). Note also that, as it was
shown in \cite[Corollary 4.4]{Ru$_6$}, a $\rho$-uniform one-sided
Markov shift $T_G$ is isomorphic to the Bernoulli shift $T_\rho$
iff the cofiltration $\xi(T_G)$ is isomorphic to standard
cofiltration $\xi(T_\rho)$.

\medskip

The purpose of this paper is to classify the $\rho$-uniform
one-sided Markov shifts.
We show that every ergodic $\rho$-uniform
Markov shift $T$ can be represented in a {\bf canonical form}
$T = T_G $ by means of a {\bf canonical} (uniquely determined by
$T$) stochastic graph $G$. In the canonical form, two such shifts
$T_{G_1}$ and $T_{G_2}$ are isomorphic if and only if their
canonical stochastic graphs $G_1$ and $G_2$ are isomorphic.

\medskip

First we consider (Section \ref{s2}) general $\rho$-uniform
endomorphisms and use the following results from \cite{Ru$_6$}.
Any ergodic $T \in \UE$ can be represented as a skew product over
$T_\rho$ on the space
$X_\rho \times Y_d$, $d \in \N \cup \{\infty\} $, where $Y_d$
consists of $d$ atoms of equal measure $\frac{1}{d}$ for
$d < \infty$ and $Y_\infty$ is a Lebesgue space with no atoms,
(see Section \ref{ss2.2} below). According to \cite{Ru$_6$} we
introduce the {\bf minimal index} $d(T)$ of $T \in \UE$ as the
minimal possible $d$ in the above skew product representation of
$T$. The index $d(T)$ is an invariant of the endomorphism $T$ and
$d(T) = 1$ iff $T$ is isomorphic to the Bernoulli shift $T_\rho$.

Other important invariants of $T \in \UE$ (introduced also in
\cite{Ru$_6$}) are the {\bf partitions} $\bldgga \bldlTr$, $\bldgb
\bldlTr$ and the {\bf index} $\bldd_{\bldgga \colon \bldgb}
\bldlTr$.
The partition $\gga (T)$ is the smallest (i.e. having
the most coarse elements) measurable partition of $X$ such that
almost all elements of the partition $\; \gb_n := \gga(T) \vee
T^{-n} \gep \;$ have homogeneous conditional measures for all $n$.
The corresponding tail partition is defined by $\; \gb (T) =
\bigwedge_{n=1}^{\infty} \gb_n  \geq \gga (T) \;$. and the index
$d_{\gga \colon \gb}(T) $ is the number of elements of $\gb(T)$ in
typical elements of the partition $\gga(T)$ (Proposition
\ref{gga:gb}).

It was proved in \cite{Ru$_6$} that $\; d(T) = d_{\gga \colon \gb}
(T) < \infty \;$ for any $\rho$-uniform one-sided Markov shift $T
= T_G$. This result implies, in particular, that $T_G$ is {\bf
simple} in the sence of Definition \ref{def SUE}.
The classification of general simple $\rho$-uniform endomorphisms
is reduced to a description of equivalent $d$-extensions of the
Bernoulli shift $T_\rho$ (Theorem \ref{simple}).

\medskip

Next we turn to $\rho$-uniform one-sided Markov shifts.

It is easy to see that a Markov shift $T_G$ is $\rho$-uniform iff
the graph $G$ satisfies the following condition: For any vertex
$u$ the set $G_u$ of all edges starting in $u$, equipped with the
corresponding weights $\; p(g) \;,\; g \in G_u \;$, is isomorphic
to $(I,\rho)$.
This means that the transition probabilities of the
Markov chain (starting from any fixed state)  coincide with $\;
\rho(i) \;,\; i \in I \;$, up to a permutation.
We call these graphs and Markov chains $\bldrho${\bf -uniform.} In
particular, $(I,\rho)$ itself is considered as a $\rho$-uniform
graph having a single vertex. The corresponding Markov shift is
the Bernoulli shift $T_\rho$.

Following \cite{AsMaTu} we use in the sequel graph homomorphisms
of the form $\; \phi \colon G_1 \raro G_2 \;$, which are assumed
to be {\bf weight preserving} and {\bf deterministic}, i.e. right
resolving in the terminology of \cite{AsMaTu}, (see Definition
\ref{hom} for details).
Thus a stochastic graph $G$ is $\rho$-uniform iff there exists a
homomorphisms $\; \phi \colon G \raro I \;$.

Two particular kinds of homomorphisms are of special interest in
our explanation, they are homomorphisms of {\bf degree 1} and {\bf
d-extensions}.
A homomorphism $\; \phi \colon G_1 \raro G_2 \;$
has degree $1$, $\; d(\phi) = 1\;$, if the corresponding factor
map $\; \Phi_\phi \colon X_{G_1} \raro  X_{G_2} \;$ is an
isomorphism.
So that $\; \Phi_\phi \circ T_{G_1} \;=\; T_{G_2}
\circ \Phi_\phi \;$, i.e. $T_{G_1}$ and $T_{G_2}$ are isomorphic.

The d-extensions homomorphism are defined in Section \ref{ss3.2}
by the condition: $\; |\phi^{-1}g| = d \;,\; g \in G \;$.
They can be described (up to equivalence) by the {\bf graph skew
products}, (see Example \ref{GSP} and Definition \ref{def GSP} in
Section \ref{ss3.2}).

\medskip

As the first step to the construction of the canonical graph we
show (Theorem \ref{phi bar}) that any homomorphism $\; \phi \colon
G \raro I \;$ can be  extended to a $d$-extension $\bphi$ by
homomorphisms of degree $1$ (See Diagram \ref{diag phi bar}).
To this end we consider a $\bldd${\bf -contractive} semigroup
$\cS(\phi)$, associated with the homomorphism $\phi$, and the
corresponding {\bf persistent} sets (Section \ref{ss4.4}). Thus we
reduce the classification problem to the study of diagrams of the
form
\begin{equation}\label{pipsi}
(\pi,\psi) \;\colon\;
\xymatrix{ \bH \ar[r]^{\pi} & H \ar[r]^{\psi} & I}
\end{equation}
where $\bH$ is a $d$-extension, $\psi$ is a degree $1$
homomorphism and the shift $T_\bH$ is isomorphic to the shift
$T_G$.

The second step is to minimize $d$ in the above Diagram
\ref{pipsi}. We show (Theorem \ref{phi bar d(T)}) that, passing
possibally to a "$n$-stringing" graph $G\pn$, one can choose the
minimal $ d = d(T)$. Note that the result is based on
\cite[Theorem 4.2 and 4.3]{Ru$_6$}.

The third final step is to reduce the homomorphism $\psi$ in
Diagram \ref{pipsi} as much as possible. Let $\Ext$ denotes the
set of all $d$-extensions of the Bernoulli graph $(I,\rho)$ of the
form (\ref{pipsi}).
We show that $\Ext$ can be equipped with a natural {\bf partial
order} "$\preceq$" and {\bf equivalence relation} "$\sim$"
(Definition \ref{partial order}).
The minimal elements of $\Ext$ with respect to the order are
called {\bf irreducible} (Definition \ref{irreduc}). We describe
these irreducible $(\pi,\psi)$-extensions by means of the
persistent $d$-partitions, associated with elements of $\Ext$
(Theorem \ref{reduc part}).

\medskip

Now we can formulate the main result of the paper
(Theorems \ref{canon form} and \ref{classification}).
\begin{itemize}
 \item {\it Let $T_G$ be a $\rho$-uniform ergodic one-sided Markov
  shift.
  A stochastic graph $\bH = \bH(G)$ is said to be a canonical graph
  for the shift $T$ if there exists an irreducible
  $(\pi,\psi)$-extension (\ref{pipsi}) from
  $\Ext$ with $\; d = d(T) \;$ such that the shift $T_\bH$ is
  isomorphic to $T_G$.}
 \item {\it Any $\rho$-uniform ergodic one-sided Markov shift can be
  represented in the canonical form $T = T_\bH$ by a canonic graph
  $\bH = \bH(G)$.}
 \item {\it In this canonical form, two shifts $T_{\bH_!}$ and
  $T_{\bH_2}$ are isomorphic iff the canonical graphs
  $\bH_1$ and $\bH_1$ are isomorphic, and iff the corresponding
  irreducible $(\pi,\psi)$-extensions are equivalent.}
\end{itemize}

\medskip

The paper is organized as follows.

In Section \ref{s2} we study general  $\rho$-uniform endomorphisms
(class $\UE$) and {\bf simple} $\rho$-uniform endomorphisms
(subclass $\SUE$).
Following \cite{Ru$_6$}, we introduce the {\bf
partitions} $\bldgga \bldlTr$, $\bldgb \bldlTr$ and the {\bf
index} $\bldd_{\bldgga \colon \bldgb} \bldlTr$.
Two main conclusions of the section are Theorem \ref{simple}
(classification of simple $\rho$-uniform endomorphisms) and
Theorem \ref{simp mar}, which states that every ergodic
$\rho$-uniform one-sided Markov shift $T_G$ is simple and $\;
d(T_G) \;=\; d_{\gga \colon \gb}(T_G) \;<\; \infty  \;$.

In Section \ref{s3} we consider general properties of stochastic
graphs and their homomorphisms.
In particular, we define $\bldrho${\bf -uniform } graphs
corresponding to $\rho$-uniform Markov shifts.
We prove that the index $d(T_G)$ of any ergodic $\rho$-uniform
Markov shift $T_G$ is finite (Theorem \ref{zet del}).
This follows from the finiteness of the degree $d(\phi)$ of any
homomorphism $\phi \colon G \raro I$ from any $\rho$-uniform graph
$G$ onto the standard Bernoulli graph $(I,\rho)$. The degree
$d(\phi)$, in turn, can be computed by means a special {\bf
d-contractive} semigroup $\cS(\phi)$, induced by $\phi$ (Theorem
\ref{d(phi)}).

Section \ref{s4} contains some essential stages of the proof of
Main Theorems \ref{canon form} and \ref{classification}.
Homomorphisms of degree $1$ and extensions of the Bernoulli graph
are considered in Sections \ref{ss4.1} and \ref{ss4.2}.
Theorem \ref{Equ ext} (Section \ref{ss4.3}) reduces the
classification of skew product over Markov shifts $T_H$ to the
classification of the corresponding graph skew product over $H$.
In Sections \ref{ss4.4} and \ref{ss4.5}, we study the set $\Ext$
of all $(\pi,\psi)$-pairs of the form (\ref{pipsi}).
The main result of Section \ref{s4} is Theorem \ref{reduc ext},
which claims the existence and uniqueness of the irreducible
$(\pi,\psi)$-pair $(\pi_*,\psi_*)$, majorized by a given
$(\pi,\psi) \in \Ext$.

In Section \ref{s5} we prove Main Theorems \ref{canon form} and
\ref{classification} and give some consequences and examples.
As a consequence we prove also (Theorem \ref{common exten 1}) that
two shifts $T_{G_1}$ and $T_{G_2}$ are isomorphic iff the graphs $G_1$ and
$G_2$ have a common extension of degree $1$.

\bigskip

We do not study here the classification problem for general, not
necessarily $\rho$-uniform, one-sided Markov shifts as well as the
classification problem of the cofiltrations, generated by the
shifts. Our approach seems to be a good tool to this end and we
hope to deal with these two problems in another paper.

We do not also consider the classification problem of one-sided
Markov shifts with infinite invariant measure, in particular, of
null-recurrent one-sided Markov shifts. One can find a good
introduction to the topic and more references in \cite[Chapters 4
and 5]{Aar}.


\bigskip

\section{ Class of $\rho$-uniform endomorphisms }
\label{s2}

\bigskip

\subsection { Lebesgue spaces and their measurable partitions }
\label{ss2.1}

We use terminology and results of the Rokhlin's theory of Lebesgue
spaces and their measurable partitions (See \cite{Rok$_1$},
\cite{Rok$_2$}). An improved and more detailed explanation can be
found in \cite{ViRuFed}. We fix the terms "homomorphism ,
isomorphism, endomorphism"  only for {\bf measure preserving} maps
of Lebesgue spaces.

Let  $(X, \cF, m)$ be a Lebesgue space with $ mX=1$.
The space $X$ is called {\bf homogeneous } if it is non-atomic or
if it consists of $d$ points of measure $\; \frac{1}{d} \;,\; d
\in \N \;$.

Let $ \gz $ be a partition of $X$ onto mutually disjoint sets $ C
\in \gz $. The element of $ \gz $  containing a point $ x $ is
denoted by $ C_\gz(x) $. The partition $ \gz $ is measurable iff
there exists a measurable function $ f \colon X \raro \R $ such
that
$$
 x \stackrel{\gz}{\sim} y \Longleftrightarrow
 C_\gz(x) = C_\gz(y) \Longleftrightarrow  f(x)=f(y) \;,\; x,y \in X
$$
Elements of $ \gz $ are considered as Lebesgue spaces
$\; (C, \cF^C , m^C) \;,\; C \in \gz \;$,
with canonical system of conditional measures
$\; m^C \;,\;C \in \gz \;$.
We shall denote also by $ m(A|C)$ the conditional measures
$\; m^C(A \cap C) \;$ of a measurable set $ A \in \cF $ in the element
$ C$ of $ \gz $.

Two measurable partitions $\gz_1$ and $\gz_2$ are said to be {\bf
independent } ($\gz_1 \perp \gz_2 $) if the corresponding
$\gs$-algebras $ \cF (\gz_1) $ and $ \cF (\gz_2)$ are independent
, where $\cF (\gz)$ denotes the $m$-completion of the
$\gs$-algebra of all measurable $\gz$-sets. We shall write also
$\; \gz_1 \perp \gz_2 \pmod \gz \;$ if the partitions $\gz_1$ and
$ \gz_2 $ are {\bf conditionally independent} with respect to the
third measurable partition $\gz$. This means that
$$
 m(A \cap B \;|\; C_\gz (x)) =
 m(A|C_\gz(x)) \cdot m(B \;|\; C_\gz(x))
$$
for all $ A \in \cF (\gz_1), B \in \cF (\gz_2) $ and a.a. $ x \in X $.

We denote by $ \gep = \gep_X $ the partition of $X$ onto separate
points and by $\nu = \nu_X$ the trivial partition of $X$.

An {\bf independent complement} of $\gd$ is a measurable partition
$ \eta $ such that $\; \gz \perp \eta \;$ and $\; \gz \vee \eta =
\gep \;$. The partition $\gz$ admits an independent complement iff
almost all elements  $ (C,m^C) $ of $\gz$ are mutually isomorphic.
The collection of all independent complements of $\gz$ is denoted
by $IC(\gz)$.

\medskip

We shall use induced endomorphisms, which are defined as follows.
Let $ A \in \cF $ , $ mA > 0$ and $T$ be an endomorphism of
$(X,m)$. Then the return function
\begin{equation}\label{ret fun}
  \gf_A(x) := min \{ n \geq 1 \colon T^{n}x \in A \}
  \;\;,\;\; x \in A
\end{equation}
is finite a.e. on $A$.
The {\bf induced endomorphism} $T_A$ on $A$ is defined now by
$\; T_Ax = T^{\gf_A(x)}x \;$.
It is an endomorphism of $ (A,\cF \cap A, m{|}_A) $ and it is ergodic
if $T$ is ergodic .

$|E|$ denotes the cardinality of the set $E$


\bigskip

\subsection {Classes ${\mathcal UE} \bldl \bldrho \bldr$ and index
$\bldd \bldlTr$.}
\label{ss2.2}

Let $(I,\rho)$ be a finite or countable state space
$$
  \rho = \{ \rho(i) \;,\; i \in I \} \;,\;\; \rho (i) > 0
  \;\;,\;\; \sum_{i \in I}{\rho(i)} = 1 .
$$

\begin{definition} \label{def rho Bern}
 An endomorphism $T$ of a Lebesgue space $(X,m)$ is said to be
 $\bldrho${\bf -unform} or {\bf finitely} $\bldrho${\bf -Bernoulli}
 endomorphism $\; (T \in \UE \;$, if there exists a discrete measurable
 partition $\gd$ of $X$, which satisfies the following condition:
 \begin{enumerate}
  \item[(i)]  $\; distr \; \gd  = \rho \;$, i.e.
   $\; \gd = \{ B(i) \}_{i \in I}$
   with $\;\;m(B(i)) = \rho (i), \;\; i \in I $,
  \item[(ii)] $\; \gd \in IC (T^{-1} \gep) $ , i.e. \ \ \
   $\gd \perp T^{-1} \gep \;$ and  $\; \gd \vee T^{-1}\gep = \gep$.
 \end{enumerate}
\end{definition}

So $ \UE $ denotes the class of all $\rho$-unform endomorphisms.
Denote by $\gD_{\rho}(T)$ the set of all partitions $\gd$
satisfying the condition $(i)$ and $(ii)$. Then $T \in \UE$ means
$\gD_\rho(T) \neq\emptyset $.

For $T \in \UE$ and $\gd \in \gD_{\rho}(T)$ define
\begin{equation}\label{delta n}
     \gd\pn = T^{-n+1} \gd \;\;,\;\;\;
     \gd\pn = \{ T^{-n+1}B(i) \}_{i \in I }  \;\;,\;\;   n \in \N
\end{equation}
Then $\; distr \;\gd_n = \rho \;$ and the partitions
$ \; \gd_1 \;,\; \gd_2 \;,\;\gd_3 \;,\; \ldots \; $ are independent.

The partitions
\begin{equation}\label{delta (n)}
  \gd\pn = \bigvee_{k=1}^{n} \gd_{k} \;\;,\;\;\;
  \gd^{(\infty)} = \bigvee_{k=1}^{\infty} \gd_k
\end{equation}
satisfy for all $n$ the conditions
$$
 \gd\pn \in IC(T^{-n} \gep) \;\;,\;\;
 \gd\pinf \perp T^{-n} \gep \pmod { \gd\pinf \wedge T^{-n} \gep }
$$
and
$$
  \gd\pinf \vee T^{-n} \gep = \gep \;\;,\;\;
  \gd\pinf \wedge T^{-n} \gep \;=\; T^{-n} \gd\pinf \;.
$$
In particular, let $ T=T_\rho$ be a Bernoulli endomorphism, which
acts on the space
$$
  (X_\rho, m_\rho) = \prod_{n=1}^{\infty} (I,\rho)
$$
as the one-sided shift
$$
  T_{\rho} x = \{ x_{n+1}\}_{n=1}^{\infty} \;\; , \;\;\;
  x = \{ x_n \}_{n=1}^{\infty} \in X_\rho \;\;.
$$
We can set
\begin{equation}\label{delta rho}
  \gd_\rho = \{ B_\rho(i) \}_{i \in I} \;\;,\;\; B_\rho(i) =
  \{ x = \{ x_{n}\}_{n=1}^{\infty} \in X_{\rho}
  \;\; \colon \;\; x_1=i \} \;.
\end{equation}
Then $\; \gd_\rho \in \gD_\rho (T_\rho )$ and $\gd_{\rho}$ is an
one-sided Bernoulli generator of $T_\rho$, that is
$$
 \gd_{\rho}^{(\infty)} = \bigvee_{n=1}^{\infty} T^{-n+1}\gd_{\rho}
  = \gep_{X_{\rho}} \;.
$$

In general case, for $T \in \UE $ and $\gd \in \gD_\rho(T)$, the
partition $\gd^{(\infty)}$ does not equal $\gep$, but we can
define the canonical factor map
$$
  \Phi_{\gd} \colon X \ni x \raro
  \Phi_{\gd} (x) = \{i_{n} (x) \}_{n=1}^{\infty} \in X_{\rho} \;,
$$
where $i_{n}(x) \in I$ is uniquely defined by the inclusion
$ T^{n}x \in B(i_{n}(x)) \in \gd $.

The homomorphism $\Phi_{\gd}$ satisfies
$\; \Phi_{\gd} \circ T = T_{\rho} \circ \Phi_{\gd} \;$
and it determines the following representation of $T$ by a skew
product over $ T_{\rho} $ (See \cite[Proposition 2.2]{Ru$_6$}).
\begin{proposition}\label{T decompos}  Let $T \in \UE $ be an
 endomorphism of $(X,m)$ and  $\gd \in \gD_\rho(T)$.
 Then
 \begin{enumerate}
  \item[(i)] There exists an independent complement $\gs$ of the
   partition $\gd\pinf$.
  \item[(ii)] The pair $\; (\gd\pinf , \gs ) \;$ induces decomposition
   of the space $(X,m)$ into the direct product
   $\; (X_\rho \times Y \;,\; m_{X_\rho} \times m_{Y})$ such that the
   factor map $\Phi_\gd$ coincides under the decomposition with the
   canonical projection
   $$
    \pi \; \colon \; X_\rho \times Y \ni (x,y) \raro x \in X_\rho
   $$
   and
   $$
    \gd = \pi^{-1} \gd_{\rho}
    \;\;, \;\; \gd\pinf = \pi^{-1} \gep_{X_\rho}
    = \gep_{X_\rho} \times \nu_Y
    \;\;,\;\; \gs = \nu_{X_\rho} \times \gep_Y
   $$
  \item[(iii)] The endomorphism $T$ is identified with the following
  skew product over $ T_\rho $
   \begin{equation}\label{skew1}
    \bT (x,y) = ( T_\rho x , A(x) y)
    \;\;,\;\; (x,y) \in X_{\rho} \times Y
   \end{equation}
   where $\; \{ A(x), x \in X_\rho \} \;$ is a measurable family of
   automorphisms of $Y$.
  \item[iv)] If $T$ is ergodic,  $Y$ is a homogeneous Lebesgue space.
 \end{enumerate}
\end{proposition}
Every homogeneous Lebesgue space $Y$ is isomorphic to $\; Y_d
\;.\;d \in \N \cup \{\infty\} \;$, where $\; Y_{\infty} \;$ is the
Lebesgue space with a continuous measure and $\; Y_d \;$,$\;d \in
\N \;$, consists of $d$ points of measure $\frac{1}{d}$. Thus for
any ergodic $T$ endomorphism $\; T \in \UE \;$ and $\; \gd \in
\gD_{\rho} (T) \;$ there exists $\; d = d(T,\gd) \in \N \cup
\{\infty\} \;$ such that
$$
 u_{\gd^{(\infty)}}(x)
 := m^{C_{\gd^{(\infty)}} (x)}(\{x\}) = \frac{1}{d}
$$
for a.a. $x \in X $.
\begin{definition}\label{def d(T)}
 \begin{enumerate}
  \item[(i)]  The number $\; d(T,\gd) \;$ will be called
   the {\bf index} of $T \in \UE $ with respect to
   $\; \gd \in \gD_{\rho}(T) \;$.
  \item[(ii)] The {\bf minimal index} $d(T)$ of $T$ is defined as
   \begin{equation}\label{d(T)}
    d(T) \;=\; min \; \{\; d(T,\gd) \;,\; \gd \in \gD_{\rho} (T)\}
   \end{equation}
 \end{enumerate}
\end{definition}
Note that an ergodic endomorphism $T$ is isomorphic to the
Bernoulli shift $T_\rho$ iff $T \in \UE$, and $d(T) = 1$, that is,
there exists $\gd \in \gD_\rho(T)$ such that $d(T,\gd) = 1$, i.e.
$\gd\pinf = \gep $.


\bigskip

\subsection{Partitions $\bldga \bldlTr$, $\bldgb \bldlTr$,
$\bldgga \bldlTr$ and indices $\bldd_{\bldga} \bldlTr$,
$\bldd_{\bldgga \colon \bldgb} \bldlTr$}
\label{ss2.3}

Let $T$ be an endomorphism of $(X,m)$ and let $ \{ \xi_n
\}_{n=1}^{\infty} $ be the decreasing sequence of measurable
partitions $\; \xi_n := T^{-n}\gep \;$, generated by $T$. The
element of $\xi_n$ , containing a point $x \in X$, has the form
$\; C_{\xi_n}(x) = T^{-n}(T^n x) \;$,

In order to introduce the partitions $\gga(T)$ and $\gb(T)$,
consider the measurable functions
$$
  u_n(x) = m^{C_{\xi_n}(x)} (C_{\xi_{n-1}}(x)) \;\;\;,
\; n \in \N \;\;\; , \; x \in X  \;,
$$
where $\xi_0 := \gep $. With these $\; u_n \colon X \raro [0,1]
\;$ we can consider the measurable partitions
$$
  \gga_n = \bigvee_{k=1}^{n}u_{k}^{-1} \gep_{[0,1]}
  \;\; ,\;\; n \in \N \;,
$$
generated by $\;u_k\;,\; k \leq n \;$,  and also
\begin{equation}\label{gga gbn gb}
 \gga = \bigvee_{n=1}^{\infty} \gga_n \;\;,\;\;
 \gb_n = \gga \vee T^{-n}\gep \;\;,\;\;
 \gb = \bigwedge_{n=1}^{\infty} \gb_n
\end{equation}
We shall write $\; \gga_n(T) \;,\; \gga(T)\;,\; \gb_n(T)\;,\;
\gb(T) \;$ to indicate $T$, if it will be necessary.
\begin{proposition} \label{gga:gb}
 Suppose that $T\in \UE $ and T is ergodic. Then there exists
 $d \in \N \cup \{\infty\}$ such that
 $$
   m^{C_{\gga}(x)}(C_{\gb}(x))=\frac{1}{d}
 $$
 for a.a. $x\in X$.
\end{proposition}

We may define now the index $\; d_{\gga \colon \gb}(T) \;$ of an
ergodic endomorphism $T \in \UE$ as the number $d$ constructed in
Proposition \ref{gga:gb}, i.e.
$$
   d_{\gga \colon \gb}(T) := (m^{C_\gga (x)}(C_\gb (x)))^{-1}
$$
for a.e. $x \in X$.

We shall use the following properties of the partitions (\ref{gga
gbn gb})
\begin{proposition}\label{gga gd d(T)}
Suppose that $\; T \in \UE \;$, let $\; \gd \in \gD_\rho(T) \;$ and
the partitions $\; \gd_n \;,\; \gd\pn \;,\; \gd\pinf \;$ defined
by (\ref{delta n}) and (\ref{delta (n)}). Then
 \begin{enumerate}
  \item[(i)] $\; \gga\pn \le \gd_n
   \;\;,\;\; \gb_n  \perp \gd\pn \pmod {\gga_n}
   \;\;,\;\; n \in \N \;.$
  \item[(ii)] $\gga \le \gd\pinf
   \;\;,\;\; \gb \perp \gd\pinf \pmod {\gga} \;.$
  \item[(iii)] $\; d_{\gga \colon \gb}(T) \; \leq \; d(T) \;$
 \end{enumerate}
\end{proposition}
We shall also use the {\bf tail} measurable partition
$\; \ga(T) := \bigwedge_{n=1}^{\infty} T^{-n} \gep \;$.
An endomorphism $T$ is called {\bf exact} if $\ga(T) = \nu$.
The {\bf tail index} $d_\ga(T)$ (which is, in fact, the {\bf
period } of $T$) is defined as follows: $d_\ga(T)= \infty$ if $
X{/}_{\ga(T)}$ is a continuous Lebesgue space and $d_{\ga}(T) = d$
if $ X {/}_{\ga((T)}$ consists of $d$ atoms of measure
$\frac{1}{d}$. So that $d_\ga(T) \in \N \cup \{\infty\}$.

It is easily to see, that
\begin{equation}\label{ga gb gga}
  T^{-1} \ga = \ga     \;\;,\;\;  \ga \vee \gga \le \gb \;\;,\;\;
  T^{-1} \gga \le \gga \;\;,\;\;  T^{-1} \gb \le \gb
\end{equation}
and $\; \ga \perp \gd\pinf \;$ for any $ \gd \in \gD_{\rho}(T) $.

Turning to the canonical projection $\Phi_\gd$ we have
\begin{proposition} \label{ga(T rho)}
 \begin{enumerate}
  \item[(i)]  $ \ga(T_\rho) = \nu
  \;\;,\;\; \gb(T_\rho) = \gga(T_{\rho}) \;$.
  \item[(ii)] $ \gga_n(T)   = \Phi_{\gd}^{-1}\gga_n(T_\rho)
   \;\;,\;\; \gga(T)   = \Phi_{\gd}^{-1}\gga(T_\rho) \;$.
 \end{enumerate}
\end{proposition}
The stated above propositions were proved in
\cite[ Propositions 2.5 - 2.9]{Ru$_6$}.


\bigskip

\subsection{Simple $\bldrho$-uniform endomorphisms}
\label{ss2.4}

We use now the partitions $\gga(T)$  and $ \gb(T)$ to introduce an
important subclass of the class $\UE$
\begin{definition} \label{def SUE}
An endomorphism $T \in \UE $ of a Lebesgue space $(X,m)$ is said
to be a {\bf simple} $\rho$-uniform endomorphism $(T \in \SUE )$,
if  there exists partition $\gd \in \gD_\rho(T)= IC (T^{-1} \gep)$
such that
 \begin{equation}\label{SUE1}
  \gd^{(\infty)} \vee \gb(T)= \gep
 \end{equation}
\end{definition}
We denote by $\SUE$ the class of all simple $\rho$-uniform
endomorphisms.

The Bernoulli endomorphism  (one-sided Bernoulli shift)
$T=T_{\rho}$  belongs to $\SUE$. In this case there exists a
partition $ \gd = \gd_{\rho} \in \gD_{\rho}(T)$ such that $
\gd\pinf = \gep$ and hence $ \gb(T) \vee \gd\pinf = \gep $
\begin{remark}\label{rem SUE} It is easily to show that the condition
(\ref{SUE1}) holds iff there exists an independent complement $\gs
\in IC(\gd^{(\infty)}) $ of $\gd^{(\infty)}$ that satisfies
 \begin{equation}\label{SUE2}
   \gs \vee \gga(T) = \gb(T)    \;\;,\;\;
   \gs \in IC (\gd^{(\infty)})  \;\;,\;\;
   \gd \in \gD_\rho(T)= IC (T^{-1} \gep)
 \end{equation}
\end{remark}
\begin{proposition}\label{simple1}
Suppose $\;T \in \UE \;$ is ergodic and $\; d(T) < \infty \;$.
Then $T$ is simple iff $\; d(T) = d_{\gga \colon \gb}(T) \;$.
\end{proposition}
\proof Since $\; d(T) < \infty \;$ we have, by Proposition
\ref{gga gd d(T)} (iii), that $\; d_{\gga \colon \gb}(T) \leq d(T)
< \infty \;$.
Definition of the index $d_{\gga \colon \gb}(T)$ (Proposition
\ref{gga:gb}) means that $\; m^{C_\gga(x)}(C_\gb(x)) = d^{-1} \;$
for a.a. $x \in X$ and $\; d = d_{\gga \colon \gb}(T) \in \N \;$.
Almost every element of $\gga(T)$ consists precisely of $d$
elements of the partition $\gb(T)$.
On the other hand there exists $\gd \in \gD_\rho(T)$ such that
almost every element of the corresponding partition $\gd\pinf $
consists precisely of $d(T)$ points, $d \leq d(T)$. By Proposition
\ref{gga gd d(T)} (ii) we have
$$
 \gb(T) \perp \gd\pinf \pmod {\gga(T)} \;\;,\;\;
 \gb(T) \wedge \gd\pinf = \gga(T) \;.
$$
Whence, the condition \ref{SUE1} holds iff $\; d(T) = d \;$.
\qed

\medskip

Let $\; T \in \SUE \;$.
By Proposition \ref{T decompos} any choice of the partition
$\gs$ in the equality (\ref{SUE2}) determines a skew product
representation (\ref{skew1}) of $T=\bT$ over $T_\rho$. Herewith,
all statements of Proposition \ref{T decompos} hold and we have
also by (\ref{SUE2}) and Proposition \ref{ga(T rho)})
\begin{equation}\label{SUE3}
 \gb(T_\rho)=\gga(T_\rho)
 \;\;,\;\; \gga(\bT) = \gga(T_\rho) \times \nu_Y
 \;\;,\;\; \gb(\bT) = \gga(T_\rho) \times \gep_Y
 \end{equation}
These arguments imply the following
\begin{theorem}\label{simple}
 Let $T$ be a $\rho$-uniform simple endomorphism, $\; T \in \SUE \;$.
 \begin{enumerate}
  \item[(i)]  $T$ can be represented in the skew product form
  (\ref{skew1}) $T=\bT$ over $T_\rho$
   \begin{equation}\label{Tbar}
    \bT (x,y) = ( T_\rho x , A(x) y)
    \;\;,\;\; (x,y) \in X_{\rho} \times Y
   \end{equation}
   with a measurable family $\; \{ A(x) \;,\; x \in X_\rho \} \;$ of
   automorphisms of $Y$ such that
   $\; \gb(\bT) = \gga(T_\rho) \times \gep_Y \;$.
  \item[(ii)] Two such skew product endomorphisms $\bT_k$, $k=1,2,$
   \begin{equation}\label{Tbar 1,2}
    \bT_k (x,y) = ( T_\rho x , A_k(x) y)
    \;\;,\;\; (x,y) \in X_{\rho} \times Y
   \end{equation}
   are isomorphic iff the corresponding families
   $\; A_1(x) \;$ and $\; A_2(x) \;$ are cohomologous, i.e.
   \begin{equation}\label{A2W=WA1}
     A_2(x) W(x) = W(T_\rho x) A_1(x) \;\;,\;\; x \in X_\rho
   \end{equation}
   for a measurable family of $\; \{ W(x) \;,\; x \in X_\rho \} \;$
   of automorphisms of $Y$.
 \end{enumerate}
\end{theorem}
\proof Part (i) follows from Proposition \ref{T decompos} with
(\ref{SUE3}).

Let $\bT_1$ and $\bT_2$ be two skew product endomorphisms of the form
(\ref{Tbar 1,2}).
Denote $\; {\ti W}(x,y):=(x,W(x)y) \;$.
Then (\ref{A2W=WA1}) implies
$ \bT_2 \circ S = S \circ \bT_1 $ if we use the automorphism
$S = {\ti W}$.

Conversely, suppose there exists an automorphism $S$ such that
$ \bT_2 \circ S = S \circ \bT_1 $. Then the partitrions
$$
 \bgga := \gga(\bT_1) = \gga(\bT_2) = \gga(T_\rho) \times \nu_Y
$$
and
$$
 \bgb := \gb(\bT_1) = \gb(\bT_2) = \gga(T_\rho) \times \gep_Y
$$
are invariant with respect to $S$. Moreover, $\bgga$ is
element-wise invariant with respect to $S$. Hence, $S {|}_C
(\bgb{|}_C) = \bgb{|}_C $ for almost every element $C \in \bgga$.
The restriction $S{|}_C$ induces a factor automorphism $W_C$ on
the factor space $ C{/}_{\bgb{|}_C} \cong Y $.
We obtain a measurable family $\; W(x) := W_{C(x)} \;,\; x \in
X_\rho  \;$, of automorphisms of $Y$. Since the partition $\;
\bgga = \gga(T_\rho) \times \nu_Y \;$ is $\bT_1)$- and
$\bT_2$-invariant, the functions $ A_1(x) $ and $ A_2(x) $ (as
well as $W(x)$) are constant on elements of $\gga(T_\rho)$.
Therefore the equality $ \bT_2 \circ S = S \circ \bT_1 $ implies $
\bT_2 \circ {\ti W} = {\ti W} \circ \bT_1 $ and (\ref{A2W=WA1})
holds. \qed

\medskip

Consider two special cases.

\medskip

\noindent {\bf Absolutely non-homogeneous $\bldrho$.} The
distribution $\; \rho = \{\rho(i) \}_{i \in I} \;$ is called
absolutely non-homogeneous if $\; \rho(i) \neq \rho(j) \;$ for all
$\;i \neq j\;$.

In this case we have $\; \gga_1(T) \vee T^{-1}\gep = \gep \;$.
On the other hand $\; \gga_1(T) \perp T^{-1}\gep \;$.
Thus $\; \gD_\rho(T) \;$ consists of the only partition, which is
$\; \gd = \gga_1(T) \;$. Hence
$$
  \gd\pinf = \gga(T) \;\;,\;\;
  \gb_n(T) = \gga_n(T) \vee T^{-n}\gep = \gep \;,
$$
$$
  \gb(T) = \bigwedge_{n=1}^{\infty} \gb_n(T) = \gep \;\;,\;\;
  \gb(T) \vee \gd\pinf = \gep
$$
Thus we have
\begin{proposition} Every $\rho$-uniform endomorphisms with absolutely
non-homogeneous $\;\rho \;$ is simple.
\end{proposition}

\medskip

\noindent {\bf Homogeneous $\bldrho$.} We have another extremal
case if $\;\rho \;$ is homogeneous, i.e. if for some $\; l \in \N
\;$,$\; I = \{1,2, \ldots ,l\;\} \;$ and $\; \rho(i) = l^{-1} \;$
for all $\; i \in I  \;.$

All the functions $ \;u_n \;$, which generate the partitions
$\; \gga_n(T) \;$, are constant,
$$
  u_n(x) \;=\; m^{C_{\xi_n }(x)} (C_{\xi_{n-1}}(x))
  \;=\; l^{-1} \;,\; n \in \N \;,\; x \in X
$$
We have $\; \gga(T) = \gga_n(T) = \nu \;,$ and $\; \gb_n(T) =
T^{-n}\gep \;$, whence, $\; \gb(T) = \bigwedge_{n=1}^{\infty}
T^{-n}\gep  = \ga(T) \;$.
Therefore, for any  $\;\gd \in \gD_{\rho}(T)\;$ the equality
(\ref{SUE1}) is equivalent to $\; \gd\pinf \vee \ga(T) = \gep \;$.
On the other hand $\; \gd\pinf \perp \ga(T) \;$ for every $\;\gd
\in \gD_{\rho}(T)\;$.

Thus we have for $ T \in \UE $ with homogeneous $\rho$

\begin{proposition} Let $ T \in \UE $ with homogeneous $\rho$. Then
 \begin{enumerate}
  \item[(i)] $T$ is simple iff there exists
   $\;\gd \in \gD_{\rho}(T)\;$ such
   that $ \gd\pinf \in IC(\ga(T)) $.
  \item[(ii)] The skew product decomposition in Theorem \ref{simple}
  is a direct product $\; T_\rho \times S \;$ with
  $ S = T {/}_{\ga(T)} $.
  \item[(iii)] Two such direct products $\; T_\rho \times S_1 \;$ and
   $\; T_\rho \times S_2 \;$ are isomorphic iff the automorphisms
   $ S_1 $ and $ S_2 $ are isomorphic.
  \item[(iv)] If, in addition, $T$ is exact, i.e. $\;\ga(T) = \nu \;$,
   then $T$ is simple iff $T$ is isomorphic to $\; T_\rho \;$.
 \end{enumerate}
\end{proposition}
It is easily to construct a skew product $T$ over $T_\rho$, which
is exact and has entropy $\; h(T) > h(T_\rho) = \log l \;$. Every
such endomorphism is $\rho$-uniform, $\; T \in \UE \;$, but it is
not isomorphic to $T_\rho$, whence, it is not simple. See also
\cite{FeR}, \cite{HeHo}, \cite{HeHoR}, \cite{Ho}, \cite{HoR}), for
more interesting examples of such kind of endomorphisms.
\begin{remark} It can be shown that there exist non-simple exact
endomorphisms in each class $\; \UE \;$ in the case, when $\rho$
is not absolutely non-homogeneous, i.e. $\; \rho(i) = \rho(j)\;$
for some $\;i,j \in I\;$.
\end{remark}

The following result plays an important role in present paper.
\begin{theorem}\label{simp mar}
Every ergodic $\rho$-uniform  one-sided Markov shift $T_G$,
corresponding to a positively recurrent Markov chain on a finite
or countable state space, is simple and
 \begin{equation}\label{d=dg:b}
  d(T_G) \;=\; d_{\gga \colon \gb}(T_G) \;<\; \infty  \;.
 \end{equation}
\end{theorem}
\proof The last statement \ref{d=dg:b} was proved in
\cite[Theorem 4.3]{Ru$_6$}.
It implies that $T_G$ is simple by Proposition \ref{simple1}.
\qed


\bigskip

\section{ Stochastic graphs and their homomorphisms. }
\label{s3}

\bigskip

\subsection {Stochastic graphs and Markov shifts}
\label{ss3.1}

We need some terminology concerning stochastic graphs and their
homomorphisms.

Consider a directed graph with countable (finite or
infinite) set $G$ of edges. Denote by $G\po0$ the set of all vertices
of the graph. We also denote by $s(g)$ the starting vertex and by
$t(g)$ the terminal vertex of an edge $g \in G$
$$
   \xymatrix{   t(g) & s(g) \ar[l]_{g} }
$$
The maps
$$
 s \; \colon \; G \ni g \; \raro \; s(g) \in G\po0 \;\; \;,\; \;\;
 t \; \colon \; G \ni g \; \raro \; t(g) \in G\po0
$$
completely determine the structure of the graph $G$,

In the sequel we assume that both the sets
$$
 _vG = \{ g \in G \;\colon\; t(g) = v \} \;\;\;,\;\;\;
 G_u = \{ g \in G \;\colon\; s(g)= u \}
$$
are not empty for all vertices $ \;u\;,v\;\in G\po0 $.

Denote by $\; G\pn \;$ the set of all paths of length $n$ in
$G$ , i.e.
\begin{equation}\label{G(n)}
  G\pn = \{ g_1 g_2 \ldots g_n \in G^n \;\colon\;
  s(g_1)=t(g_2) , \ldots , s(g_{n-1})=t(g_n) \}
\end{equation}
A graph $G$ is said to be {\bf irreducible} if for every pair of
vertices $ \;u\;,v\;\in G\po0 $ there exists a finite  $G$-path
$ \;g_1g_2 \ldots g_n \in G\pn \;$ such that $\; u=s(g_n) \;$ and
$\; v=t(g_1) \;$.

Take into account that we use here and in the sequel the notation
 $\; g_1 \; g_2\; \ldots \; g_n \;$ for {\bf backward} paths
$$
 \xymatrix{
    t(g_1)  &  s(g_1) = t(g_2) \ar[l]_-{g_1}
            &  s(g_2) = t(g_3) \ar[l]_-{g_2}
            &                  \ar[l]_-{g_3}    }
               \;\;\dots\;\;
 \xymatrix{ &  s(g_n)          \ar[l]_-{g_n}    }
$$
A graph $G$ is called {\bf stochastic} if its edges $g$ are equipped
with positive numbers $ p(g) $ such that $\; \sum_{g \in
G_u}\;p(g) \;=\; 1 \;$ for all $\; u \; \in G\po0 \; $. The
weights $\; p(g) \;,\; g \in G \;$, determine the backward
transition probabilities of the Markov chain induced by $G$.

We shall assume in the sequel that there exist stationary
probabilities $ p\po0(u) > 0 $ on $ G\po0$ such that
\begin{equation}\label{PR}
 \sum_{u \in G\po0}\;p\po0(u) \;=\; 1 \;\;\;,\;\;\;
 \sum_{g \in _vG} p(g) p\po0 (s(g)) \;=\; p\po0 (v)
\end{equation}
for all vertices $ \;u\;,v\;\in G\po0 $.

It is known that the stationary probabilities on $G\po0$ exist iff
the corresponding to $G$ Markov chain is positively recurrent.
If, in addition, the irreducibility condition hold, the stationary
probabilities $\;p\po0(u) \;,\; u \in G\po0 \;$ on the vertices
are uniquely determined by the transition probabilities $\;
p(g)\;,\;g \in G\;$ on the edges.

Thus any stochastic graph $(G,p)$ induces a Markov chain on the
state space $G$ with the transition probabilities matrix $\; P
\;=\;{(P(g,h))}_{g\in G , h \in G} \;$, where
$$
P(g,h) \;=\;
\left\{
 \begin{array}{ll}
  p(h) \;, & if \;\; t(g)=s(h) \; \\
  0  \;\;, & otherwise.
 \end{array}
\right.
$$
In the sequel we mainly deal with stochastic graphs, which induce
irreducible positively recurrent Markov chains.

\medskip

The one-sided Markov shift $ T_G $, induced by the stochastic
graph $ G $ , is defined as follows. Let
$$
 X_G = \{x={\{g_{n}\}}_{n=1}^{\infty} \in G^{\N} \;\colon\;
 s(g_1)=t(g_2)\;,\; s(g_2)=t(g_3) \;,\;  \dots \; \}
$$
 and the Markov measure $m_G$ on $ X_G $ is given by
$$
 p(g_{1}) \; p(g_{2}) \; \ldots \; p(g_n) \; p\po0(s(g_n))
$$
on the cylindric sets of the form
$$
 A(g_1\;g_2\; \ldots \;g_n)
 \;:=\; \{x={\{x_k\}}_{k=1}^{\infty} \in X_G
 \;\colon\; x_1=g_1 \;,\; \ldots \;,\; x_n=g_n \;\}
$$
where $\; g_1\;g_2\; \ldots \;g_n) \; \in \; G\pn \;$ is a
$G$-path of length $n$ in $G$.

The one-sided shift $T_G$ acts on the probability space
$\;( X_G , m_G) \; $  by
$$
 T_G({\{x_n\}}_{n=1}^{\infty}\;) \;=\; {\{x_{n+1}\}}_{n=1}^{\infty}
$$
and $T_G$ preserves the Markov measure $m_G$.
The shift $T_G$ is ergodic iff the graph $G$ is irreducible. Under
the irreducibility condition, the stationary probabilities $p\po0$
on $G\po0$ and, hence, the $T_G$-invariant Markov measure $m_G$
are uniquely determined by the stochastic graph $\; (G,p) \;$.

The coordinate functions
$$
 Z_n \;\colon\; X_G \ni x={\{x_k\}}_{k=1}^{\infty}\; \raro \; x_n
 \in G \;\;,\;\;n \in \N
$$
form a stationary Markov chain on $ ( X_G , m_G) $ with the
backward transition probabilities
$$
P(g,h) \;=\;m_G\{\; Z_n = h \;|\; Z_{n+1} = g \;\} \;=\; p(h)
\;\;,\;\; n \in \N
$$
for all $\; (h,g) \in G^{(2)} \;$.

\medskip

Consider now the partitions
$$
 \gz_n = {Z_n}^{-1}\gep_G =
 T_G^{-n+1}\gz_1 = \{ T_G^{-n+1}A(g) \}_{g \in G}
 \;\;,\;\; n \in \N \;,
$$
generated by  $Z_n$ on $ X_G$, where
$$
 A(g) \;=\;
 \{\; x={\{x_k\}}_{k=1}^{\infty} \in X_G \;\colon\; x_1=g \;\}
$$
Setting $ \gz = \gz_1 $ and $ T = T_G $, we have
\begin{equation}\label{gen par}
 \bigvee^{\infty}_{n=1}T^{-n+1}\gz = \gep
\end{equation}
\begin{equation}\label{Mar par}
  \gz  \perp  \bigvee^{\infty}_{n=1}T^{-n}\gz \pmod { T^{-1}\gz }
\end{equation}

Recall that a measurable partition $\gz$ of $(X,m)$ is said to be
a {\bf one-sided Markov generator} or {\bf one-sided Markov
generating partition} for an endomorphism T of $(X,m)$, if the
above conditions (\ref{gen par}) and (\ref{Mar par}) hold.

The partition $ \gz_G $ will be called the {\bf standard}
one-sided Markov generator of the one-sided Markov shift $T_G$ on
$ X_G $.
\begin{example}[{\bf Standard Bernoulli Graph}]\label{Bern graph}
 Let $(I,\rho)$ be a finite or countable alphabet and
 $$
  \rho = \{ \rho(i) \;,\; i \in I \} \;,\;\; \rho (i) > 0
  \;\;,\;\; \sum_{i \in I}{\rho(i)} = 1 .
 $$
 be a probability on $I$.
 We shall consider $(I,\rho)$ as a stochastic graph, which has the set of
 edges $i \in I$ with weights $\; \rho(i) \;$ and a single vertex, denoted
 by$\;"o"\;$.
 So $G\po0 = \{o\}$ is a singleton and $ s(i) = t(i) = o $ for all
 $i \in I$.
 We shall say that $(I,\rho)$ is the {\bf standard Bernoulli graph.}

 \medskip

 For instance,
 $\; \xymatrix{ o \ar@(ul,dl) []_{p} \ar@(dr,ur) []_{q} } \;$
 if $\; |I| = 2 \;$ and $\; \rho = (p,q) \;$.

 \medskip

 The corresponding to $(I,\rho)$ one-sided  Markov shift $T_I$ coincides
 with the Bernoulli shift $ T_I = T_\rho $.
 The generating partition $ \gz_I$ coincides with the standard Bernoulli
 generator $\; \gd_\rho = \{ B_\rho(i) \}_{i \in I} \;$, defined by
 (\ref{delta rho}).
\end{example}

\medskip

\noindent
{\bf Induced shift} $\; \bldT_\bldu \;$.
For any $\; u \in G\po0 \;$, denote
$$
 D(u) \;:=\; \{\; x={\{x_k\}}_{k=1}^{\infty} \in
 X_G \;\colon\; t(x_1)=u \;\} \;\;,\;\; u \in G\po0 \;.
$$
and consider the partition
$\; \gz\po0 := {\{ D(u) \}}_{u \in G\po0} \;$
on the space $X_G$.
The partition $\gz\po0$ is a Markov partition with respect to shift $T_G$ , i.e.
\begin{equation}\label{Mar par 0}
  \gz\po0  \perp  T_G^{-1} \gep_{X_G} \pmod { T_G^{-1}\gz\po0 } \;,
\end{equation}
but it is not a one-sided generator for $T_G$, in general.

 We shall use in the sequel the endomorphisms $ T_u := (T_G)_{D(u)} $,
induced by the shift $T_G$ on elements $ D(u) $ of $\gz\po0$,
$\; u \in G\po0 \;$.
The Markov property (\ref{Mar par 0}) provides that for every
$u \in G\po0 $ the induced endomorphism $ T_u $ is a Bernoulli
shift. More exactly, in accordance with the general definition of
return functions (\ref{ret fun}) we have
$$
 \gf_u(x)=\gf_{D(u)}(x) := min \{ n \geq 1 \colon T_G^{n}x \in D(u) \}
 \;\;,\;\; x \in D(u)
$$
and
$$
  T_u x \;=\; T_G^{\gf_u(x)} x \;\;,\;\; x \in X_G \;.
$$
Take $\; I_u = \bigcup_{n=1}^\infty I_{u,n} \;$, where $I_{u,n}$ be the
set of all $ g_1g_2 \ldots g_n \in G\pn $ such that
\begin{equation} \label{I u n}
 t(g_1) = s(g_n) = u \;\;,\;\;
 s(g_k) = t(g_{k+1}) \neq u \;,\; k = 1,2, \ldots ,n-1 \;.
\end{equation}
Define also $\; \rho_u = \{\rho_u(i)\}_{i \in I_u} \;$ by
\begin{equation} \label{rho u}
\rho_u(i) = p(g_1)  p(g_2) \ldots  p(g_n) \;\;,\;\;
        i = g_1g_2 \ldots g_n \in I_{u,n} \;.\; n \in \N \;.
\end{equation}
For any $\; i = g_1g_2 \ldots g_n \in I_{u,n} \;$ we set
$ B_u(i) := A(g_1g_2 \dots g_n)$ and consider the partition
$\; \gz_u = \{ B_u(i)\}_{i \in I_u} \;$, whose elements are enumerated
by the alphabet $I_u$.
The Markov property (\ref{Mar par 0}) implies that the partitions
$\; T_u^{-n} \gz_u \;,\; n \in \N \;$
are independent. Thus
\begin{proposition} \label{Tu}
The induced endomorphism $T_u$ is isomorphic to the Bernoulli
shift $T_{\rho_u}$ and $ \gz_u $ is a one-sided Bernoulli
generator of$T_u$.
\end{proposition}


\bigskip

\subsection {Graph homomorphisms and skew products}
\label{ss3.2}

Now we want to establish the class of graph homomorphisms that we shall use.
\begin{definition} \label{hom}
Let $\; G \;$ and $\; H \;$ be two stochastic graphs.
 \begin{enumerate}
  \item[(i)] A map $\; \phi \colon G \raro H \;$ is a {\bf graph
   homomorphism} if there exists a map
   $\; \phi\po0 : G\po0 \raro H\po0 \;$ such that
   $$
    s(\phi (g)) = \phi\po0 (s(g)) \;\;\;,\;\;\;
    t(\phi (g)) = \phi\po0 (t(g))
   $$
   for all $\; g \in G \;$.
   (Note that, if the map $\; \phi\po0 \;$ exists it is unique.)
  \item[(ii)] A graph homomorphism $\;\phi \colon G \raro H \;$ is
   {\bf deterministic} if $\; \phi\po0 (G\po0) = H\po0 \;$
   and for every $\; u \in G\po0 \;$ the restriction of $\phi $ on
   $ G_u$
   $$
    \phi {|}_{G_u} \;\colon\; G_{u} \; \rightarrow \;H_{\phi\po0 (u)}
   $$
   is a bijection of this set onto $\; H_{\phi\po0 (u)} \;$.
  \item[(iii)] A graph homomorphism is {\bf weight preserving} or
   $\bldp${\bf-preserving} if  $\; p(\phi (g)) \;=\; p(g) \;$ for all
   $\; g \in G \;$.
 \end{enumerate}
\end{definition}
Two edges $g_1$ and $g_2$ are said to be {\bf congruent}
if $$
 s(g_1)=s(g_2) \;\;,\;\; t(g_1)=t(g_2) \;\;,\;\; p(g_1)=p(g_2) \;.
$$
The map $ \phi\po0 $ in the above definition is uniquely
determined by $ \phi $, but $ \phi\po0 $ does not determines $
\phi$ if $G$ has congruent edges.

Anyway one can use a more explicit notation
$$
 ( \phi , \phi\po0 ) \;\colon\; (G , G\po0) \; \raro \; (H , H\po0 )
$$
for the homomorphism $\;\phi \colon G \raro H \;$.

We shall denote by $\; \Hom (G,H) \;$ the set of all weight
preserving deterministic graph) homomorphisms $\;\phi \colon G
\raro H \;$. In the sequel the term "homomorphism" always means
just {\bf weight preserving deterministic graph homomorphism}.
\begin{proposition}\label{factor map}
Let $\; \phi \colon G \raro  H\;$ be a map.
 \begin{enumerate}
  \item[(i)] If $ \phi $ is a graph homomorphism, it induces a factor
  map
   $$
   \Phi_\phi \;\colon\; X_G  \raro  X_H \;\;,\;\;
   \Phi_\phi ( {\{x_n\}}_{n=1}^{\infty} ) \;=\;
   {\{ \phi ( x_n ) \}}_{n=1}^{\infty}
   $$
   such that $\; \Phi_\phi \circ T_G \;=\; T_H \circ {\Phi}_{\phi}\;.$
  \item[(ii)] If, in addition, $\phi $ is weight preserving, the factor map
   $ \Phi_\phi $ is measure preserving,
   $\; ( m_H = m_G \circ {\Phi_\phi}^{-1} ) $.
  \item[(iii)] If $\;\phi\;$ is also deterministic, the shift $T_G$
  can be represented as a skew product
   \begin{equation}\label{skew prod}
   \bT (x,y) = ( T_H x \;,\; A(x) y) \;\;,\;\; (x,y) \in X_H \times Y
   \end{equation}
   where $\; \{ A(x) \;,\; x \in X_H \} \;$ is a measurable family of
   automorphisms of $Y$.
  \item[(iv)] If $T_G$ is ergodic,  $Y$ is a homogeneous Lebesgue
  space.
 \end{enumerate}
\end{proposition}
\proof Parts (i) and (ii) follow directly from Definition
\ref{hom}. Part (iii) and (iv) can be proved by analogy with
Proposition \ref{T decompos} \qed

\medskip

Moreover
\begin{theorem}\label{fin ind}
 Let $\phi \in  \Hom (G,H) $ and suppose that the shift
 $ T_G $ is ergodic.
 Then there exists $\;d \in \N \;$ such that
 $\; | \Phi_\phi^{-1}(x) | \;=\; d \;$ for almost all
 $\; x \in X_H \;$.
 That is, in the skew product (\ref{skew prod}) the space $Y$ is
 finite, $\; |Y| = d \;$.
\end{theorem}
Note that Theorem \ref{fin ind} claims the finiteness of $d$ even
in the case, when the graph $G$ is not finite, i.e. $ |G| = \infty
$. This is a consequence of positive recurrence of the
corresponding to $G$ Markov chain. The skew product decomposition
(\ref{skew prod}) of $T_G$ over $T_H$ is a $d$-extension.

Theorem \ref{fin ind} was proved earlier in a particular case,
when $H$ is a Bernoulli graph, i.e. when $ H\po0 = \{ o \} $ is a
singleton (See Theorem 3.3 and Corollary 3.4 from \cite{Ru$_6$},
and also Theorem \ref{zet del} below).

We omit the proof of Theorem \ref{fin ind} in general case , since
only the pointed out particular case is considered in this paper.
\begin{definition}\label{def d(phi)}.
The integer $d$ in Theorem \ref{fin ind}, i.e. the degree of the
factor map $\Phi_\phi$ will be called the {\bf degree} of the
homomorphism $\phi$.
\end{definition}
Denoting the degree by $\; d(\phi) \;$, we have
$\; d(\phi) = |\Phi_\phi^{-1}(x)| \;$ for a.a. $ x \in X_H $.

\medskip

The following construction plays a central role in our explanation.
\begin{example}[{\bf Graph Skew Product}]\label{GSP}

Let $\;d \in \N \;$ and let $\; Y_d = \{ 1,2, \ldots , d \} \;$
consists of $d$ points of measure $\frac{1}{d}$. Denote by $\; \Ad
= {\mathcal A}(Y_d) \;$ the full group of all permutations of
$Y_d$.

Given a stochastic graph $ H $, equipped with a function
$\; a \colon H \ni h \raro a(h) \in \Ad \;$,
we construct a stochastic graph $ \bH_a $ and a homomorphism
$\; \pi_H \colon \bH_a  \raro  H \;$ by
$$
  \bH_a \;=\; H \;\times \; Y_d  \;\; , \;\;
  {\bH_a}\po0 \;=\; H\po0 \; \times \; Y_d \; ,
$$
with
$$
  s(\bh) = (s(h),y)     \;\;,\;\;
  t(\bh) = (t(h),a(h)y) \;\;,\;\;  p(\bh) = p(h)
$$
for $\; \bh = (h,y) \in \bH_a = H \times Y_d \;$ and also
$$
  p\po0(\bu) = p\po0(u) \;\;,\;\;
  \bu = (u,y) \in {\bH_a}\po0 = H\po0 \times Y_d
$$.
The natural projection
$$
 \pi_H \colon \bH_a = H \times Y_d \raro H \;\;,\;\;
 {\pi_H}\po0 \colon \bH_a^{(0)} = H\po0 \times Y_d \raro H\po0
$$
is a homomorphism.
\begin{definition}\label{def GSP}
We shall say that the graph $ \bH_a $ is a {\bf skew product} over
$H$ and the homomorphism $\;\pi_H \colon \bH_a \raro H \;$ is a {\bf
graph skew product} (or {\bf GSP})  $d$-extension of $H$ .
\end{definition}
\end{example}
In the above construction we have $\; | {\pi_H}^{-1} (h) | = d \;$
for all $\; h \in H \;$
and this is, in fact, a characteristic property of the graph
skew product $d$-extension in the following sense
\begin{definition}\label{hom equi}
 Two homomorphisms $\; \phi_k \colon G_k \raro H \;,\; k = 1,2 \;$ are
 said to be {\bf equivalent} if $\; \phi_2 = \gk \circ \phi_1 \;$ for
an appropriate isomorphism $\; \gk \colon G_1 \raro G_2 \;$.
\end{definition}
\begin{definition}\label{d-uniform} Let $\; d \in \N \;$.
 A homomorphism $\; \phi \in \Hom (G,H) \;$ is called a
 {\bf $d$-extension} if
 \begin{equation}\label{phi-1=d}
  |{\phi}^{-1} (h)| = d \;\;,\;\; h \in H \;.
 \end{equation}
\end{definition}
\begin{proposition}\label{d-unif=GSP}
 Any $d$-extension $ \phi \colon G \raro H $ is equivalent to a
 GSP $d$-extension $\; \pi_H \colon \bH_a \raro H \;$.
\end{proposition}
\proof Let $\; \phi \in \Hom (G,H) \;$ is a $d$-extension.
Since $\phi$ is deterministic the restrictions $\; \phi{|}_{G_u}
\;$ are bijections between $G_u$ and $ H_{\phi\po0(u)} $ for all $
u \in H\po0 $. Hence the condition (\ref{phi-1=d}) is equivalent
to
$$
 |{\phi\po0}^{-1}(u)| = d \;,\; u \in H\po0 \;.
$$
For each $u \in H\po0$ we can choose a bijection $ w_u $ of
$ {\phi\po0}^{-1} (u) $ onto $ Y_d $.
With any fixed choice of these bijections we set
$$
 H \ni h \raro  a(h) = w_{t(h)} \circ {w_{s(h)}}^{-1} \in \Ad \;,
$$
and consider the corresponding skew product graph $ \bH_a $.
The bijections $ w_u $ uniquely determine an isomorphism
$\; \gk \colon G \raro \bH_a \;$ such that
$\; \phi =\pi_H \circ \gk \;$.
\qed

\begin{remark}\label{re GSP shift} The Markov shift $T_{\bH_a}$
corresponding to a graph skew product $ \bH_a $ can be identified with
the skew product endomorphism $ \bT_{H,a} $, defined by
$$
  \bT_{H,a} (x,y) = (T_H x , {a(x_1)}^{-1}y) \;\;,\;\;
  x ={\{ x_n \}}_{n=1}^{\infty} \in X_H \;,\; y \in Y_d \;.
$$
and thus any $d$-extension is a homomorphism of degree $d$.

Indeed, the shift $T_{\bH_a}$ acts on the space
$$
 X_{\bH_a} = \{ {\{ (x_n,y_n) \}}_{n=1}^{\infty} \;\colon\;
 x = {\{ x_n \}}_{n=1}^{\infty} \in X_H \;,\;
 y_n = a(x_n)y_{n+1} \in Y_d \}
$$
and the map
$$
  \Psi \;\colon\; X_{\bH_a} \ni {\{ x_n \}}_{n=1}^{\infty} \raro
  ({\{ x_n \}}_{n=1}^{\infty},y_1)  \in X_H \times Y_d
$$
realizes the identification, that is,
$\; m_H \otimes m_{Y_d} = m_{\bH_a} \circ \Psi^{-1} \;$ and
$\; \bT_{H,a} \circ \Psi = \Psi \circ T_{\bH_a} \;$.
Note also that
\begin{equation}\label{Psi gz}
 \Psi \; \gz_{\bH_a} \;=\; \gz_H \times \gep_{Y_d} \;\;,\;\;
 \Psi \; \gz_{\bH_a}\po0  \;=\; \gz_H
\po0 \times \gep_{Y_d} \;.
\end{equation}
\end{remark}

\medskip

Consider now two skew product endomorphisms $\; \bT_{H,a_k} \;$,
corresponding to graph skew products $ \bH_{a_k} $ with two functions
$\; a_k \colon H \raro \Ad \;,\; k=1,2 \;$.
\begin{definition}\label{cohom}
 \begin{enumerate}
  \item[(i)] Two functions $\; a_k \colon H \raro \Ad \;,\; k=1,2 \;$, are
   said to be {\bf cohomologous} with respect to $H$ if there exists
   a map $\; w \colon H\po0 \raro \Ad \;$ such that
    \begin{equation}\label{a cohom}
     a_2(h) w(s(h)) = w(t(h)) a_1(h) \;,\; h \in H
    \end{equation}
  \item[(ii)] Two measurable functions
   $\; A_k \colon X_H \raro \Ad \;,\; k=1,2 \;$ are said to be {\bf
   cohomologous} with respect to $T_H$ if there exists a measurable
   map $\; W \colon X_H \raro \Ad \;$ such that
   \begin{equation}\label{A cohom}
    A_2(x) W(x) = W(T_H x) A_1(x) \;,\; x \in X_H
   \end{equation}
 \end{enumerate}
\end{definition}
In accordance with Definitions \ref{hom equi} and \ref{cohom} we
can say now that the homomorphisms
$\;\pi_H \colon \bH_{a_k} \raro H $ are equivalent iff the
functions $\; a_k : H \raro \Ad \;,\; k=1,2 \;$ are cohomologous
with respect to $H$.

The equality (\ref{a cohom}) is equivalent to (\ref{A cohom}) if we
take
$$
 A_k(x) = {a_k(x_1)}^{-1} \;,\; k=1,2 \;\;,\;\; W(x) = w(t(x_1))
$$
for $\; x = {\{ x_n \}}_{n=1}^{\infty} \in X_H \;$ and given $a_k$
and $w$. Hence if $a_1$ and $a_2$ are cohomologous with respect to
$H$, then $A_1$ and $A_2$ cohomologous with respect to $T_H$.

We shall show in Section \ref{ss4.3} that the inverse is also
true.
\begin{remark}\label{triv exten}
Let $\chi \colon H \raro H_1 $ be a homomorphism and
$ \pi_1 : \bH_1 \raro H_1 $ be a $d$-extension of $H_1$ generated
by a function $a_1 : H_1 \raro \Ad $. Setting $ a(h) :=
a_1(\chi(h)) $ we obtain a $d$-extension $ \pi \colon \bH := \bH_a
\raro H $ of $H$. The map $ \bchi(h,y) := (\gk(h),y) \;,\; (h,y)
\in \bH$ is a homomorphism and the diagram
\begin{equation}\label{diag ext}
   \xymatrix{
    \bH \ar[d]^{\pi} \ar[r]^{\bchi} & \bH_1 \ar[d]^{\pi_1} \\
        H \ar[r]^{\chi}                    &  H_1 }
\end{equation}
commutes.
The homomorphism $ \bchi $ is called a {\bf trivial} extension
of $\bchi$. If, in addition, $ d(\chi) = 1 $, then $ d(\bchi) = 1 $
and hence the corresponding endomorphisms $\bT_{H,a}$ and
$\bT_{H_1,a_1}$ are isomorphic.
\end{remark}


\bigskip

\subsection{Stochastic $\bldrho$-unform  graphs}
\label{ss3.3}

We continue to consider $\;(I,\rho)\;$ as the standard Bernoulli
stochastic graph, (Example \ref{Bern graph})
\begin{definition} \label{rho uni}
A stochastic graph $ (G,p) $ is called $\bldrho${\bf-uniform} if
there exists a homomorphism  $\; \phi \in \Hom (G,I) \;$.
\end{definition}
For any such homomorphism  $\phi$ and for every $ u \in G\po0 $
$$
  \phi \;{|}_{G_u} \;\colon\;
  ( G_u \;,\; p \;{|}_{G_u} )\; \raro \; (I,\rho)
$$
is a weight preserving bijection.
Thus the atomic probability spaces $\; ( G_u , p {|}_{G_u} ) \;$
are isomorphic to $(I,\rho)$ for every $ u \in G\po0 $.
\begin{proposition} \label{p3.4}
 $\; T_G \in \UE \;$ iff $G$ is $\rho$-uniform.
\end{proposition}

\proof Consider the partition $\; \xi_1 := T_{G}^{-1} {\gep}_{X_G} \;$
generated by the shift  $T_G$.
The Markov property of the measure $m_G$ on $X_G$ implies
$$
 m^{C_{\xi_1}(x)}(\{x\}) \;=\; m_G \{\; Z_1=x_1 \;|\; Z_2=x_2 \;\}
  \;=\; p(x_1)
$$
for a.a. $\; x={\{x_{n}\}}_{n=1}^{\infty} \in X_G \;$,
Here $\; m^{C_{\xi_1}(x)}(\{x\}) \;$ is the conditional measure of
the point $x$ in the element $C_{\xi_1}(x) = T_G^{-1}T_Gx $ of the
partition $ \xi_1 $.
Hence for every $ u \in G\po0 $ almost all elements $(C,m_C)$ of
the partition $ \xi_1 $ are isomorphic to $\; ( G_u , p {|}_{G_u})
\;$ on the set $\; \{ x={\{x_n\}}_{n=1}^{\infty} \in X_G
\;\colon\; s(x_1) = u \} \;$. But $\;T_G \in \UE \;$ iff a.a.
elements $(C,m_C)$ are isomorphic to $(I,\rho)$. Hence $\;T_G \in
\UE \;$ iff $\;( G_u , p{|}_{G_u})\;$ are isomorphic to $(I,\rho)$
for every $ u \in G\po0 $. \qed

\medskip

Let $G$ be a $\rho$-uniform graph and $\; \phi \in \Hom(G,I) \;$.
Consider the partition
$\; \phi^{-1} \gep_I = \{ \phi^{-1} (i) \;,\; i \in I \;\}$ of $G$.
The first coordinate function
$$
 Z_1 \;\colon\; X_G \ni x={\{x_k\}}_{k=1}^{\infty}
 \; \raro \; x_1 \in G
$$
generates the following partition
$$
  \gd_\phi \;=\; Z_1^{-1} ( \phi^{-1} \gep_I )
$$
of the space $ X_G  $.
Elements of $\; {\gd}_{\phi} \;$ have the form
$$
   B(i) = Z_{1}^{-1} ( \phi^{-1} (i) ) =
   \{ x={\{x_{k}\}}_{k=1}^{\infty} \in X_G  \;\colon\;
   {\phi}(g) = i \} \;,\; i \in I
$$
Using the standard Markov generator
$$
 \gz_G = Z_{1}^{-1}\gep_G = {\{A(g)\}}_{g \in G} \;\;,\;\;
 A(g) = Z_{1}^{-1}(g)
$$
of $T_G$, we have
$$
  B(i) = {\bigcup}_{g \in {\phi}^{-1}(i)} A(g)
$$
and
$$
 m_G(B(i)) = {\sum}_{g \in {\phi}^{-1}(i)} p(g)  p\po0(s(g))
           = \rho(i) {\sum}_{u \in G\po0} p\po0(u) = \rho(i)
$$
for $\; i \in I \;$.
Hence for $\; \gd \;=\; \gd_\phi \;$ we have
\begin{equation}\label{del phi}
 \gd \in IC (\; T_{G}^{-1} \gep_{X_G} \;) =
 \gD_\rho( T_G )  \;\;,\;\;   \gd \leq \gz_G
\end{equation}
Denoting by $\; \gD_\rho(T_G , \gz_G) \;$ the set of
all $ \gd $ that satisfy (\ref{del phi}), we have also
\begin{proposition}\label{d(T,zeta)}
 $\; \gD_\rho(T_G , \gz_G) \;$ is precisely the set of all
 $ \gd $ of the form $\; \gd \;=\; \gd_\phi \;$.
\end{proposition}

\medskip

Now we introduce a semigroup $\; \cS (\phi) \;$ of maps
$\; f \colon G\po0 \raro G\po0 \;$ induced by the homomorphism $\phi$.

Let $\; i \in I \;$.
Since $\; \phi \;$ is deterministic the restriction
$\; \phi {|}_{G_u} \;\colon\; G_u \; \raro I \;$ is a bijection of
$G_u$ onto $I$ for every $u \in G\po0$. Hence for any pair $(i,u)$
there exists an unique $\; g_{i,u} \;$ such that $\; \phi (
g_{i,u} ) = i \;$ and $\; s( g_{u,i} ) = u \;$. Putting $\; f_iu =
g_{i,u} \;$, we get a map $\; f_i \colon G\po0 \raro G\po0 \;$.
Let  $\; \cS (\phi) \;$ be the semigroup generated by the maps $\;
\{ f_i \;,\; i \in I \} \;$.

Let $ \cF \cS (I) $ be the set of all finite words $\; i_1i_2
\dots i_n \;$  in the alphabet $I$.
We shall consider $ \cF \cS (I) $ as a free semigroup with the
generating set $I$ and with juxtaposition multiplication:
$$
 i_1i_2 \dots i_m \cdot j_1j_2 \dots j_n
 \;=\; i_1i_2 \dots i_mj_1j_2 \dots j_n
$$
and set
$$
 f_{i_1i_2 \dots i_n} \;=\; f_{i_1} \circ f_{i_2} \circ \; \dots \;
  \circ f_{i_{n}} \;\;,\;\; i_1i_2 \dots i_n \; \in \; I^n
$$
Then $\; i_1i_2 \dots i_n \raro  f_{i_1i_2 \dots i_n} \;$ is a
homomorphism from the semigroup $ \cF \cS (I) $ onto the semigroup
$$
\cS (\phi) = \{ f_{i_1i_2 \dots i_n} \;,\;i_1i_2 \dots i_n
\in \cF \cS (I) \} \;,
$$
generated by $\; \{ f_i \;,\; i \in I \} \;$.

\medskip

Now we can describe the partitions
$$
 \gd_\phi = \{ B(i) \}_{i \in I} \;\;,\;\;
 \gd_\phi\pn \;=\; \bigvee_{k=1}^{n} T_{G}^{-k+1} \gd_\phi
 \;,\; n \in \N
$$
as follows.

First recall that the partition $\; \gz\po0$ consists of the atoms
$D(u) = Z^{-1}(_uG) \;,\;u \in G\po0 \;$ and rename the elements
$\; A(g) \;,\; g \in G \;$ of the partition $\gz_G $ by
$$
 D(i,u) := A(g_{i,u}) \;\;,\;\; u \in G\po0 \;,\; i \in I \;.
$$
Then for all $i \in I$ and $u \in G\po0$ we have
$\; D(i,u) = B(i) \cap T_{G}^{-1}D(u) \;$,
$$
 D(i,u) = \{ x={ \{ x_n \} }_{n=1}^{\infty} \in X_G
 \;\colon\; t(x_1) = u \;,\; \phi(x_1) = i \}
$$
and
\begin{equation}\label{Bi Diu}
 B(i) = \bigcup_{u \in G\po0} D(i,u) \;\;\;,\;\;\;
 D(u) = \bigcup_{v \colon f_i(v)=u} D(i,v) \;.
\end{equation}
Further for any $\; g_1g_2 \dots g_n  \in  G\pn \;$ there exists a
unique pair $\; (i_1i_2 \dots i_n,u) \in I^n \times G\po0 \;$ such
that
\begin{equation}\label{u i t}
 u = s(g_n) \;,\; i_k = \phi (g_k) \;,\; t(g_k) = f_{i_k} (s(g_k))
 \;,\; k = 1,2, \dots ,n
\end{equation}
Hence any atom $\; A(g_1g_2\; \dots \,g_n) \;$ of the partition
$\; \gz_G\pn = \bigvee_{k=1}^{n} T^{-k+1} \gz_G \;$ can be renamed
by $\;D(i_1i_2 \dots i_n,u) = A(g_1g_2\; \dots \,g_{n}) \;$, where
the pair $\; (i_1i_2 \dots i_n,u) \;$ satisfies (\ref{u i t}).
By (\ref{Bi Diu}) any atom $\; B(i_1i_2\; \dots \;i_n) \;$ of the
partition $\; \gd_\phi\pn = \bigvee_{k=1}^{n} T_{G}^{-k+1}
\gd_\phi \;$ has the form
$$
 B(i_1i_2\; \dots \;i_n)
 \;=\; \bigcup_{u \in G\po0} D(i_1i_2 \dots i_n,u)
$$
and since
$$
 D(i_1i_2...i_n,u) \;=\; B(i_1i_2\; \dots \;i_n) \cap T^{-n} D(u)
$$
we have
$$
 m_G(B(i_1i_2\; \dots \;i_n)) \;=\;
 \rho (i_1) \rho (i_2) \; \dots \; \rho (i_n),
$$
$$
 m_G( D(i_1i_2 \dots i_n,u) ) \;=\;
 \rho (i_1) \rho (i_2) \; \dots \; \rho (i_{n})p\po0(u).
$$
Any $\gz_G\po0$-set has the form
$$
 D(E) = \{ x={ \{ x_n \} }_{n=1}^{\infty} \in X_G
 \;\colon\; t( x_1 ) \in E \}.
$$
for a subset $\; E \subseteq G\po0 \;$.
Then for any $ i_1i_2 \; \dots \; i_n \in I^n $
$$
 D(E) \cap B(i_1i_2\; \dots \;i_n) \;=\;
 \bigcup_{u \colon f_{i_1i_2 \;\dots\; i_n} (u) \in E }
 D(i_1i_2 \dots i_n,u).
$$
Hence
\begin{equation}\label{D(E)}
 m_G ( D(E) \; | \; B(i_1i_2 \; \dots \; i_n) ) \;=\;
 p\po0 ( f^{-1}_{i_1i_2\; \dots \; i_n} (E) ).
\end{equation}
Next theorem is basic for our explanation. Let
$$
 \gd_{\phi}^{(\infty)} \;=\; \bigvee_{n=1}^{\infty} \gd_\phi\pn
 \;=\; \bigvee_{n=1}^{\infty} T_{G}^{-n+1} \gd_\phi
$$
\begin{theorem} \label{zet del}
 Let $\; \phi \in \Hom (G,I) \;$. Then
 \begin{enumerate}
  \item[(i)] $\; \gz_G \; \vee \; \gd_{\phi}^{(\infty)}
  \;=\; \gep_{X_G} \;.$
  \item[(ii)] $\; d(T_G) \;\leq\; d(T_G,\gd_{\phi}) \;=\;
   d(\phi) \;<\; \infty \;$
 \end{enumerate}
\end{theorem}
\proof  It was proved in \cite[Theorem 3.3]{Ru$_6$} that if $\gz$
is a one-sided Markov generator of $T$ and
$$
 T \in  \UE \;,\; \gd \in \gD_\rho(T) \;,\; \gd \; \leq \; \gz \; ,
$$
then $\; \gz \vee  \gd^{(\infty)} = \gep \;$. Hence (i) follows by
putting $ \gd = \gd_\phi $ and $ \gz = \gz_G $ .

Since the partition $ \gz_G $ is finite or countable the equality
(i) implies that almost all elements $(C,m_C)$ of the partition
$\gd_{\phi}^{(\infty)}$ are atomic.
Taking in to account the ergodicity of $T_G$, we see that almost
all elements of $\gd_{\phi}^{(\infty)}$ consist of $d$ atoms of
measure $\frac{1}{d}$ for an natural $d$. Herewith by Definitions
\ref{def d(T)} and \ref{def d(phi)} we have $\; d =
d(T_G,\gd_{\phi}) = d(\phi)$ and, whence, (ii) follows. \qed

\medskip

We need the following sharp version of Part (i) of Theorem
\ref{zet del}
\begin{lemma}\label{zet0 del}
 $\;  \gz\po0_G \; \vee \; \gd_{\phi}^{(\infty)} \;=\; \gep_{X_G} \;$
\end{lemma}
\proof
Choose an increasing sequence of positive numbers $c_n > 0$ and an
increasing sequence of finite subsets $E_n$ of $G\po0$ such that
 \begin{equation}\label{zet0 del1}
  \bigcup_{n=1}^\infty E_n = G\po0 \;\;\;.\;\;\;
  \sum_{n=1}^\infty(1-c_n) < \infty \;\;\;.\;\;\; p\po0(E_n)>c_n \;.
 \end{equation}
Since $|E_n|<\infty$ there exist
$\; i\pn_1i\pn_2 \dots i\pn_{k_n} \in I^{k_n} \;$ and
$\; f_n := f_{i\pn_1i\pn_2 \dots i\pn_{k_n}} \in \cS (\phi) \;$
such that
$$
 |f_n(E_n)| \;=\; \min \{ |f(E_n)| \;\colon\; f \in \cS(\phi) \}.
$$
The choice of $f_n$ provides that all restrictions
$\; f{|}_{f_n(E_n)} \;,\; f \in \cS(\phi) \;$ are bijections.

Consider the sets
\[
\begin{split}
 & B_n := B(i\pn_1i\pn_2 \dots i\pn_{k_n}) \;\;,\;\; \\
 & B'_n := B_n \bigcap T_G^{-k_n}D(E_n) \;=\;
   \bigcup_{u \in E_n}D(i\pn_1i\pn_2 \dots i\pn_{k_n},u)
\end{split}
\]
and also
$$
 F_n \;:=\; \{ x \in X_G \;\colon\;
 T^{\gw_n(x)+n}x \in B^\prime _n  \} \;,
$$
where
$$
\gw_n(x) \;:=\; \min \{k \geq 0 \;\colon\; T^{n+k}x \in B_n  \}
$$
Then it is not hard to see that
$$
 m_G(F_n) \;=\; m_G(B'_n \;|\; B_n) \;=\; p\po0(E_n) \;>\; c_n \;.
$$
Set $\; F \;:=\; \liminf_{n \raro \infty}F_n \;$.
Then we have $\; m_G(F)=1 \;$, since $\;\sum(1-c_n) < \infty \;$.
By constructing, the set $F$ has the following property.
Suppose $\; x={\{x_k\}}_{k=1}^{\infty} \;$ and
$\;y={\{y_k\}}_{k=1}^{\infty} \;$ belong to $F$ and
$$
 \; \Phi_{\gd_\phi}(x) \;=\; \Phi_{\gd_\phi}(y)
 \;=\; {\{i_k\}}_{k=1}^{\infty}
 \;\in \; X_\rho \;.
$$
If $\; s(x_m)  \neq s(y_m) \;$ for some $m \geq 1 $, then
$$
 t(x_k) = f_{i_k}s(x_k) \neq t(y_k) = f_{i_k}s(y_k) \;\;,\;\;
 k \;=\; 1,\,2,\;\dots \;m \;.
$$
In other words, if $t(x_1) = t(y_1)$ and
$\; \Phi_{\gd_\phi}(x) \;=\; \Phi_{\gd_\phi}(y) \;$ then $x=y$.
Thus $\; \gz\po0_G \; \vee \; \gd_{\phi}^{(\infty)} \;=\;
\gep_{X_G} \;$ on the set $F$ of measure $1$. \qed


\bigskip

\subsection{Semigroup $ \cS (\bldphi)$ and persistent $\bldd$-sets.}
\label{ss3.4}

Let $U$ be a finite or countable set.
\begin{definition}\label{L set}
Let $ \cS$ be a semigroup of maps $\;f \;\colon\; U \;\raro \; U\;$ on
$U$ and let $\; d \in \N \;$.
Call the semigroup $\cS$ $\bldd${\bf -contractive}
if there exists a subset $\; L\; \subseteq U \;$ such that
\begin{enumerate}
 \item[(i)] $ |f(L)| \;=\; |L| \;=\; d \;$ for all $\; f \in \cS \;.$
 \item[(ii)] For every finite subset $\; E \subset U \;$ there exists
  $\; f \in \cS \;$ with $\; f(E) \subseteq L \;$.
\end{enumerate}
The sets $L$, satisfying (i) and (ii), will be called {\bf persistent
$\bldd$-sets} with respect to $ \cS $.
\end{definition}
Denote by $\; \cL (\cS) \; $ the set of all such $L$.
We have directly from the definition:
\begin{itemize}
 \item For $\; L \in \cL (\cS ) \;$ and $\; f \in \cS \;$ the
  restriction $\; f {|}_{L} \;\colon\; L \;\raro \; f(L) \;$ is a
  bijection and $\; f(L) \; \in \; \cL ( \cS ) \;$.
 \item The semigroup $ \cS $ acts transitively on $ \cL ( \cS ) $,
  i.e. for every pair $\; L_1 \;,\;L_2\; \in \; \cL (\cS ) \;$ there
  exists $\; f \in \cS \;$ such that $\; f(L_1) \;=\; L_2 \;$.
 \item The integer $d$ is equal to
  \begin{equation}\label{d(G)}
    d( \cS ) \; := \; \sup_{E \subseteq U \;\colon\; |E| < \infty}
   \;\;  \min_{ f \in \cS } \;\; |f(E)| \; .
  \end{equation}
  and $\; d(\cS) \;=\; \min_{ f \in \cS } \;\; |f(U)| \;$ if
  $\; |U| \;<\; \infty \;$.
\end{itemize}

\medskip

Let $G$ be a $\rho$-uniform stochastic graph and $ \phi \in  \Hom
(G,I) $ be a homomorphism $ \phi \colon G \raro I $.

Return to the semigroup $\; \cS(\phi) \;$ which acts on
$\; U = G\po0 \;$ .
\begin{theorem}\label{d(phi)}
Let $T_G$ be an ergodic one-sided Markov shift corresponding to a
$\rho$-uniform stochastic graph $G$ and let
$ \phi \in \Hom (G,I)$. Then the semigroup $ \cS(\phi)$ is
$d$-contractive on $G\po0 $ and
 \begin{equation}\label{d=d(G)}
  d \;=\; d(\cS(\phi)) \;=\; d(T_G , \gd_\phi) \;=\; d(\phi)
 \end{equation}
\end{theorem}
\proof To prove the theorem we shall use the partition $\gz\po0_G$
on $G\po0$.
Recall that $\gz\po0_G$ consists of all atoms of the
form $\; D(u) = Z_1^{-1} (_uG) \;,\; u \in G\po0. \;$. For any
subset $E$ of $G\po0 $ we denote
$$
D(E) = \{ x={\{x_{n}\}}_{n=1}^{\infty} \in X_G  \;\colon\;
 t( x_1 ) \in E \} = \bigcup_{u \in E} D(u) \;,
$$
i.e. $ D(E) $ is a $ \gz\po0_G$-set corresponding to $E$ in
the space $X_G$.

It follows from Theorem \ref{zet del} Part (ii) that almost all
elements $(C,m_C)$ of the partition $\gd_{\phi}^{(\infty)}$ are
isomorphic to $Y_d$, where $\; d = d(T_G,\gd_{\phi}) \in \N \;$.
Hence
$$
 m( \{x\} \;|\; C_{\gd_{\phi}^{(\infty)}}(x)) \;=\; \frac{1}{d}
$$
for a.e. $x \in X_G$.
Then Lemma \ref{zet0 del} implies that there
exists a measurable family $\; \{ l(x) \;,\; x \in X \} \;$ of
subsets $\; l(x) \subseteq G\po0 \;$ such that
\begin{equation}\label{l(x)=d}
 m ( D(l(x)) \;|\; C_{\gd_{\phi}^{(\infty)}} (x) ) \;=\; 1
 \;\;,\;\;  |l(x)| \;=\; d
\end{equation}
almost everywhere on $ X_G $.

For any $ L \subseteq G\po0 $ denote
$$
 \check{L} \;:=\; \{ x \in X_G \;\colon\; l(x) = L \} \;\;,\;\;
 \cL \;:=\; \{ L \subseteq G\po0 \;:\; m_G ( \check{L} ) > 0 \} \;,
$$
i.e. $\cL$ is the (finite or countable) set of all essential values
of the function $\;x \;\raro \; l(x) \;$. We show that $\; \cL
\subseteq \cL (\cS (\phi)) \;$, i.e. that every $\; L \in \cL \;$
satisfies the conditions (i) and (ii) of Definition \ref{L set}.

Take any finite subset $E \subseteq G\po0$ and choose $\; c > 0
\;$ such that
$ \; c  \;<\; \min_{u \in E}p\po0(u) \;$. For $\; L \in \cL \;$
and almost all $\; x={\{x_{n}\}}_{n=1}^{\infty} \in \check{L} \;$
we have by (\ref{D(E)})
$$
 lim_{n \raro \infty} m(D(L) \;|\; C_{\gd_\phi\pn}(x) \}
 \;=\; m_G (D(l(x)) \;|\; C_{\gd_\phi\pinf} (x)\} \;=\; 1
$$
and by (\ref{D(E)})
$$
 m(D(L) \;|\; C_{\gd_\phi\pn} (x)\} \;=\;
 p\po0( f^{-1}_{x_1x_2\; \ldots \;x_n} (L)).
$$
Hence we can choose $n$ and $\; (x_1x_2\; \ldots \;x_{n}) \in G\pn
\;$ such that
$$
 m(B(x_1x_2\; \dots \;x_n) \cap \check{L} ) > 0
$$
and then
$$
 p\po0 (f^{-1}_{x_1x_2\; \dots \;x_{n}} (L)) \;>\; 1-c
$$
The choice of $c$ provides
$\; f^{-1}_{x_1x_2\; \ldots \;x_n}(E) \supseteq L \;$
and thus Part (ii) of Definition \ref{L set} holds. Part (i)
follows from the equalities
\begin{equation}\label{f l(x)}
 f_{x_1x_2\; \dots \;x_{n}} l(T_G^nx) \;=\; l(x) \;\;,\;\; |l(x)| = d
 \;\;,\;\; x \in X_G
\end{equation}
by the definition of $l(x)$ . We have proved the inclusion $\; \cL
\subseteq \cL (\cS (\phi)) \;$ , which implies that the semigroup
$\cS$ is $d$-contractive with $\; d = d(T_G , \gd_\phi) $ \qed

\begin{corollary}\label{L=L(phi)}
 $\; \cL = \cL (\cS (\phi)) \;$,
\end{corollary}
\proof It was proved above that
$\; \cL \subseteq \cL (\cS (\phi)) \;$.
Take $ M \in \cL (\cS (\phi)) $ and $ L \in \cL $.
Since also $ L \in \cL (\cS (\phi))$, there exists $i_1i_2 \dots
i_{n} \in I^n $ such that
$$
\; f_{x_1x_2\;...\;x_{n}} (M) \;=\; L
\;=\; l(x) \;,\; x \in \check{L} \;
$$
Then the equality (\ref{f l(x)}) implies $\; M = l(T^{n}x) \;$ on
a set of positive measure in $X$ and hence $\; M \in \cL \;$.
Thus $\; \cL \supseteq \cL (\cS (\phi)) \;$.
\qed

\begin{remark}
Note that the notion of $d$-contractive semigroup was introduced
in \cite {Ru$_6$}, where an analog of Theorem \ref{d(phi)} was
also proved. Definition \ref{L set} is a generalization of what is
called "point collapsing" by M. Rosenblatt \cite{Ro$_1$} ,
\cite{Ro$_2$} in the case $ |U| < \infty $ and $ d=1 $ .
The case , when $ |U| = \infty $ and $ d=1 $, was considered in
\cite{KuMuTo}.
\end{remark}


\medskip

\subsection{ Graph skew product representation.}
\label{ss3.5}

From now on let $G$ be a $\rho$-uniform stochastic graph, which is
irreducible and satisfies the positive recurrence condition.
\begin{theorem}\label{phi bar}
Let $\; \phi \in \Hom (G,I) \;$ be a homomorphism of degree
$\; d = d(\phi) \;$.
Then there exists a commutative diagram
\begin{equation}\label{diag phi bar}
 \xymatrix{  \bH \ar[d]^{\bphi}  \ar[r]^{\bpsi} & G \ar[d]^{\phi} \\
                                H \ar[r]^{\psi} & I }
\end{equation}
where the graph $ \bH = \bH_a $ is a graph skew product over $H$,
generated by a function $\; a \colon H \ni h \raro a(h) \in \Ad
\;$, the homomorphism $ \bphi $ coincides with the natural
projection $\pi_H$, and both the homomorphisms $\; \bpsi \in \Hom
(\bH,G) \;$ and $\; \psi \in \Hom (H,I)\;$ are of degree 1. In
particular, $\; (\bphi,\psi) \in \Ext \;$
\end{theorem}
\proof We construct a commutative diagram
\begin{equation}\label{diag phi hat}
 \xymatrix{
  \hH \ar[d]^{\hphi} \ar[r]^{\hpsi} & G \ar[d]^{\phi} \\
                    H \ar[r]^{\psi} & I }
\end{equation}
such that the homomorphism $\; \hphi \in \Hom (\hH,H) \;$ is a
$d$-extension (See Definition \ref{d-uniform}) and both
homomorphisms $\; \hpsi \in \Hom (\hH,G) \;$ and $\; \psi \in \Hom
(H,I) \;$ are of degree 1 .

We shall use the persistent $d$-sets  $ \cL = \cL (\cS (\phi)) $
of the semigroup $ \cS (\phi) $, described in Theorem \ref{d(phi)}
(Section\ref{ss3.4}).
Since $ \cL $ is finite or countable we can enumerate the set by
an alphabet $J$, setting  $ \cL = \{ L_j \;,\; j \in J \} $.

Recall that the semigroup $ \cS (\phi) $ is $d$-contractive with
$ d = d(\cS (\phi)) $ (Theorem \ref{d(phi)}).
For any pair $\; i \in I \;,\; j \in J \;$ we have $\; |f_i(L_j)|
= |L_j| =d \;$ and the restrictions $\; f_i {|}_{L_j} \;$ is a
bijection of $\; L_j \;$ onto $\;f_i(L_j) \;$, whence, $\;f_i(L_j)
\in \cL \;$ for all $i$ and $j$. For any $\;i \in I \;$ denote by
$f_i^J$ the map $\; J \raro J \;$, which is defined by $\; f_i^J j
= j' \;$, where $\; f_i(L_j) = L_{j'} \;$.

To construct Diagram \ref{diag phi hat} define first
$\; \psi \colon H \raro I \;$
with $ \; H := I \times J \;$,$\; H\po0 := J \;$, where
$$
 s(i,j) = j \;\;,\;\; t(i,j) = f_i^J j \;\;,\;\; p(i,j) = \rho(i)
 \;\;,\;\; \psi (i,j) = i \;.
$$
Next set
$$
 \hH\po0 \;:=\;
 \{ (j,u) \in J \times G\po0 \;\colon\; j \in J \;,\; u \in L_j \}
 \;\;,\;\; \hH \;:=\; I \times \hH\po0
$$
with
$$
 s(i,j,u) = (j,u) \;\;,\;\; t(i,j,u) = (f_i^J j , f_i u) \;\;,\;\;
 p(i,j,u) = \rho(i) \;.
$$
Finally, we define the maps $\hpsi$ and $\hphi$  by
$$
 \hphi \colon \hH \ni (i,j,u) \raro (i,j) \in H \;\;,\;\;
 \hpsi \colon \hH \ni (i,j,u) \raro g_{i,u} \in G \;.
$$
where $g_{i,u}$ is uniquely determined by the conditions $s(g) =
u$ and $\phi(g) = i$.

It follows directly from this constructing that $H$ and $\hH$ are
stochastic graphs, and that $\hphi$ , $\hpsi$ and $\psi$ are
homomorphisms, and that Diagram \ref{diag phi hat} commutes. Point
out only that $\hphi$ is a $d$-extension, since $ |L_j| = d $ for
all $j$ and hence $\hphi$ is of degree $d$.
This implies that $\hat{\psi}$ and $\psi$ are of degree 1 , since
$\phi$ is of degree $d$.

It remains to apply Proposition \ref{d-unif=GSP} to the  homomorphism
$\; \hphi \colon \hH \raro H \;$.
\qed

\medskip

Next we construct d-extensions with a minimal possible $d$. Let
$G$ as above be a $\rho$-uniform stochastic graph.
Recall that $G\pn$ denotes the set of all $n$-paths in $G$, see
(\ref{G(n)}). We shall consider $G\pn$ as a stochastic graph with
the set of vertices $G^{(n-1)}$, where for any $\; g\pn =
g_1\;g_2\; \ldots \;g_n \;$
$$
 s(g\pn) \;=\; g_2\;g_3\; \ldots \;g_n \;\;,\;\;
 t(g\pn) \;=\; g_1\;g_2\; \ldots \;g_{n-1}
$$
and $\; p(g\pn) \;=\; p(g_1) p(g_2) \ldots p(g_n) \;$. If $G$ is
$\rho$-uniform, the "$n$-stringing"
graph $G\pn$ is also $\rho$-uniform.
The natural projection
$$
 \pi\pn \colon G\pn \ni g\pn =(g_1\;g_2\;
 \ldots \;g_n) \raro g_1 \in G
$$
is a homomorphism and $\; \phi \circ \pi\pn \in \Hom (G\pn,I) \;$
for any $\; \phi \in \Hom (G,I) \;$. However, if  $(I,\rho)$ has
congruent edges there exist $\; \phi_1 \in \Hom (G\pn,I) \;$,
which are not of the above form $\; \phi_1 = \phi \circ \pi\pn
\;$. It is an obvious fact, that $d(\pi\pn) = 1$, i.e. $\;
\Phi_{\pi\pn} : X_{G\pn} \raro X_G$ is an isomorphism. We use the
index $d(T,\gd)$ and the minimal index $d(T)$, which were defined
by Definition \ref{def d(T)}.

\begin{theorem}\label{phi bar d(T)}
Let $G$ be a $\rho$-uniform stochastic graph, which is irreducible
and satisfies the positive recurrence condition.
Then there exist an integer $n \in \N$, a homomorphism $\; \phi
\in \Hom (G\pn,I) \;$ and a commutative diagram
\begin{equation}\label{diag bar H Gn}
 \xymatrix{
 \bH \ar[d]^{\bphi} \ar[r]^{\bpsi} & G\pn \ar[d]^{\phi} \\
                               H \ar[r]^{\psi} & I   }
\end{equation}
such that
 \begin{enumerate}
  \item[(i)] The graph $\bH = \bH_a $ is a skew product over a graph
   $H$, generated by a function
   $\; a \colon H \ni h \raro a(h) \in \Ad \;\;,\;\;d \in \N \;$,
   and the homomorphism $ \bphi $ coincides with the natural
   projection $\pi_H$ of $\bH$ onto $H$,
  \item[(ii)] $\; d = d(\phi) = d(T_G) \;$,
  \item[(iii)] The homomorphisms $\; \bpsi \in \Hom (\bH,G) \;$ and
   $\; \psi \in \Hom (H,I) \;$ are of degree 1 .
 \end{enumerate}
\end{theorem}
\proof Let $\gz = \gz_G $ be the standard Markov generator of the
shift $T_G$. It was proved in \cite[Theorem 4.2 and 4.3]{Ru$_6$}
that there exist  $ n \in \N $ and $ \gd \in \gD_\rho(T_G) $ such
that
$$
 \gd \leq  \gz\pn := \bigvee_{k=1}^{n} T^{-k+1} \gz
 \;\;,\;\; d(T,\gd) = d(T) \;.
$$
Take $\gz\pn$ and $G\pn$ instead $\gz$ and $G$ in Proposition
\ref{d(T,zeta)}. Then we obtain $\; \phi \in \Hom (G\pn,I) \;$
with $\; \gd = \gd_\phi \;$ and thus, by using Corollary \ref{phi
bar}, we complete the proof. \qed


\bigskip

\section{Homomorphisms and finite extensions.}
\label{s4}

\bigskip

\subsection {Homomorphisms of degree 1}
\label{ss4.1}

Let $H$ be a $\rho$-uniform graph and consider a homomorphism
$\; \psi \colon H \raro I \;$.
Suppose that $\psi$ is of degree 1.
By Theorem \ref{d(phi)} the semigroup $\cS (\psi)$, generated by
$\; f_i = f_i^\psi \;,\; i \in I \;$, is $1$-contractive and all its
persistent sets are singletons.
Using $\psi$ we can identify the graph $H$ with $I \times J$, where
$\; J = H\po0 \;$ and for any $ h = (i,j) \in H $
$$
 \psi (h)=i \;,\; s(h)=j \;\;,\;\; t(h)= f_i j
 \;\;,\;\; p(h)=\rho (i)
$$
Since $d(\psi)=1$ the partition $\gd_\psi$ is a one-sided
Bernoulli generator for the Markov shift $T_H$.
The factor map $\; \Phi_{\gd_\psi} \;\colon\; X_H \raro X_\rho \;$
is an isomorphism, $\; \Phi_{\gd_\psi} \;\circ\; T_H = T_{\rho}
\;\circ\; \Phi_{\gd_\psi}$  and we can consider the Markov
partitions $\; \gz_\rho := \Phi_{\gd_\psi}(\gz_H) \;$ and $\;
\gz\po0_\rho := \Phi_{\gd_\psi}(\gz_H\po0) \;$ for $T_\rho$ on
$X_\rho$, which correspond to the Markov partitions $\gz_H$ and
$\gz_H\po0$ for $T_H$ on $X_H$.
The partition $\gd_\rho = \Phi_{\gd_\psi}(\gd_H) $ coincides with
the standard Bernoulli generator of the Bernoulli shift $T_\rho$.

Thus we have, with the notations from Section \ref{ss3.3},
$$
 \gd_\rho = \{ B_\rho (i) \}_{i \in I} \;\;,\;\;
 \gz_\rho = \{ D_\rho (i,j) \}_{ (i,j) \in I\times J} \;\;,\;\;
 \gz\po0_\rho = \{ D_\rho (j) \}_{ j \in J}
$$
where
$$
 D_\rho (i,f_i(j)) = B_\rho (i) \cap T_\rho^{-1} D_\rho (j)
 \;\;,\;\;  (i,j) \in I\times J
$$
Hence the homomorphism $\psi$ is determined by the partitions $\;
\gd_\rho \;,\; \gz_\rho \;,\; \gz\po0_\rho \;$ uniquely up to
equivalence (see Definition \ref{hom equi}).

\medskip

Our aim now is to construct a {\bf common extension of degree 1}
for two homomorphisms of degree 1.
\begin{theorem}\label{degree 1} Let
$\;  \psi_1 \colon H_1 \raro I  \;,\; \psi_2 \colon H_2 \raro I
\;$ be two homomorphisms of $\rho$-uniform graphs $H_1$ and $H'_2$
onto the Bernoulli graph $ (I,\rho) $ and suppose that $ \psi_1 $
and $ \psi_2 $ are of degree 1.
Then there exist a $\rho$-uniform graph $H$ and homomorphisms
$\psi$, $\chi_1$ and $\chi_2 $ of degree 1, for which the
following diagram commutes:
\begin{equation}\label{diag degree 1}
 \xymatrix{
  & H_2  \ar[rr]^{\psi_2} & &  I \\
  H \ar[ur]^{\chi_2} \ar[urrr]^{\psi} \ar[rr]_{\chi_1} & &
  H_1 \ar[ur]_{\psi_1} &  }
 \end{equation}
\end{theorem}
The homomorphism $\psi$ will be called a {\bf common extension} of
$\psi_1$ and $\psi_2$ of degree $1$. \proof Denote by $ (\gz_1 ,
\gz\po0_1) $ and $ (\gz_2 , \gz\po0_2) $ the pairs of Markov
partitions of the space $X_\rho$, which correspond to the
homomorphisms $\psi_1$ and $\psi_2$. Here we omit the subscript
"$\rho$" and mark the partitions and their elements by subscripts
"$1$" and "$2$", respectively.

We have to construct the desired $H$ and $\; \psi \colon H \raro I
\;$ by means of the partitions
$$
 \gz := \gz_1 \vee \gz_2 \;\;,\;\;
 \gz\po0 := \gz_1^{(0)} \vee \gz_2^{(0)} \;.
$$
By the identification $H_1= I \times J_1$ and $H_2= I \times J_2$,
we have
$$
 \gz_1 = \{ D_1 (i,j_1) \}_{ (i,j_1) \in I \times J_1} \;\;,\;\;
 \gz_1^{(0)} = \{ D_1 (j_1) \}_{ j_1 \in J_1}  \;\;,
$$
$$
 \gz_2 = \{ D_2 (i,j_2) \}_{ (i,j_2) \in I \times J_2} \;\;,\;\;
 \gz_2^{(0)} = \{ D_2 (j_2) \}_{ j_2 \in J_2}  \;\;,
$$
and then the partition $ \gz\po0 $ consists of all elements
$$
 D(j) = D_1 (j_1) \cap D_2 (j_2) \;\;,\;\; j=(j_1,j_2) \in J \;.
$$
where the set $J$ is defined by
\begin{equation}\label{eq J}
 J := \{ j=(j_1,j_2) \;\colon\; p\po0 (j) :=
 m_\rho ( D_1 (j_1) \cap D_2(j_2)) >0 \} \subset J_1 \times J_2 \;.
\end{equation}
For any $i \in I$ and $ j=(j_1,j_2) \in J $ we set
$ f_i j := (f_{1,i}  j_1 ,f_{1,i} j_2 )$.
Then
\[
\begin{split}
   D(f_i j) :=
 & D_1 (f_{1,i} j_1 ) \cap D _2(f_{2,i} j_2)
     \supseteq  D_1(i,j_1) \cap D_2(i,j_2)  =      \\
 & B(i) \cap T_\rho^{-1} (D_1 ( j_1) \cap D_2( j_2)) =
   B(i) \cap T_\rho^{-1} D(j)
\end{split}
\]
Since $ \gd $ and $ T_\rho^{-1} \gep $ are independent, this implies
$$
  p\po0 (f_i j) = m_\rho ( D(f_i j)) \geq
 m_\rho(B(i) \cap T_\rho^{-1} D(j)) = \rho(i) p\po0 (j) \;.
$$
Hence $f_i j \in J$ for all $j \in J$ and $ i \in I $ .

Thus we are able to define a stochastic graph $H := I \times J $ with
$H\po0 := J$ such that for any $j \in H\po0 $ and $h=(i,j) \in H $
$$
 s(h) := j \;,\; t(h)
 := f_i (j) \;,\; p(h) := \rho (i) \;,\; \psi (h) := i \;.
$$
The construction provides that $H$ is a $\rho$-uniform graph ,
$p\po0$ is a stationary probability on $H\po0 $ and $\psi \colon H
\raro I $ is a homomorphism of index $1$. Moreover, if we set
$$
  \chi_1(h) := (i,j_1) \;,\; \chi_2(h) := (i,j_2) \;\;,\;\;
  h=(i,j_1,j_2) \in H= I \times J_1 \times J_2 \;,
$$
then $\chi_1 \colon H \raro H_1$ and $\chi_2 : H \raro H_2$ are
homomorphisms and Diagram \ref{diag degree 1} commutes. \qed

\medskip

We shall use also the following sharpening of the previous
theorem, which can be proved in a similar way.
\begin{theorem}\label{degree 1 sharp} Let
$$
  \gk_1 \colon H_1 \raro H_0  \;\;,\;\; \gk_2 \colon H_2 \raro H_0
   \;\;,\;\; \psi_0 \colon H_0 \raro I
$$
be homomorphisms of $\rho$-uniform graphs $H_1 \;,\; H_2$ and
$H_0$ and suppose they are of degree 1.
Then there exist a $\rho$-uniform graph $H$ and homomorphisms
$\chi$, $\chi_1$ and $\chi_2 $ of degree 1, for which the
following diagram commutes
\begin{equation}\label{diag degree 1 sharp}
 \xymatrix{
 & H_2  \ar[rr]^{\gk_2} & & H_0 \ar[rr]^{\psi_0} & & I \\
 H \ar[ur]^{\chi_2} \ar[urrr]^{\chi} \ar[rr]_{\chi_1} & &
   H_1 \ar[ur]_{\gk_1} & & & }
 \end{equation}
\end{theorem}
Note that this theorem holds without adding of homomorphism
$\psi_0$ i.e. for graphs, which are not necessary $\rho$-uniform,
but we do not use the fact in this paper.


\bigskip

\subsection{ Extensions of Bernoulli graphs.}
\label{ss4.2}

Consider a very special case of the graph skew product
construction $\bH_a$ (see Example \ref{GSP}), when the graph $H$
is the standard Bernoulli graph $(I,\rho)$. Let $ d \in \N $ and
let $\; a \colon I \raro \Ad \;$ be a function on $I$ with the
values $\; a(i) , i \in I , \;$ in the group $\Ad$ of all
permutations of $\; Y_d = \{1,2, \ldots ,d \} \;$.
According to the general GSP-construction we have $\; \bI_a = I
\times  Y_d \;$ , $\; {\bI_a}\po0 = Y_d \;$ and $\; \pi \colon
\bI_a \raro I \;$, where for any $ \bh = (i,y) \in \bI_a $
$$
 s(\bh) = y \;,\; t(\bh) = a(i)y \;,\; \pi (\bh) = i \;,\;
 p(\bh) = \rho(i) \;,\, p\po0 (y) = d^{-1} \;.
$$
The stochastic graph $\bI_a$ is $\rho$-uniform and it is
irreducible iff the group $\gG (a) $, generated by $\; \{ a(i) , i
\in I \} \subseteq \Ad \;$, is transitive on $Y_d$.

As it was noted in Section \ref{ss3.2} (see Remark \ref{re GSP shift}) the
Markov shift $T_{\bI_a} $ is isomorphic to the skew product
$\bT_{\rho,a}$ , which acts on $X_\rho \times Y_d$ by
\begin{equation}\label{bT rho a}
  \bT_{\rho,a} (x,y) = (T_\rho x , {a(x_1)}^{-1} y) \;\;,\;\;
  x ={\{ x_n \}}_{n=1}^{\infty} \in X_\rho \;,\; y \in Y_d \;.
\end{equation}
\begin{theorem}\label{th rho cohom} Let
$\; \pi_k \;\colon\; \bI_{a_k} \raro  I \;,\; k = 1,2, \;$ be two
$d$-extensions of the Bernoulli graph $(I,\rho)$, generated by
functions $ a_k \colon I \raro \Ad $, respectively. Let the
functions $\; A_k \colon X_\rho \raro \Ad \;,\; k=1,2, \;$ are
defined by
\begin{equation}\label{Am(x)}
 A_k(x) := {a_k(x_1)}^{-1}  \;\;,\;\; x
 ={\{x_n\}}_{n=1}^{\infty} \in X_\rho \;.
\end{equation}
If there exists a measurable function $ W \colon X_\rho \raro \Ad
$ such that
\begin{equation}\label{eq rho cohom}
 A_2(x) \cdot W(x) = W(T_\rho x) \cdot A_1(x) \;\;,\;\; x \in X_\rho
\end{equation}
then $W(x)$ does not depend on $x$ , i.e. $W(x) = w_0 \in \Ad $
a.e. on $X_\rho$. Thus $A_1$ and $A_2$ are cohomologous with
respect to $T_\rho$ iff $a_1$ and $a_2$ are conjugate in $\Ad$,
i.e. $\; a_2(i) \cdot w_0 = w_0 \cdot a_1(i) \;,\; i \in I \;$.
\end{theorem}
Note that the last equality means the equivalence of $a_1$ and
$a_2$ in the sense of Definition \ref{cohom} , since $I\po0 = \{ o
\}$.

\medskip

To proof the theorem we need the following simple lemma.
\begin{lemma}\label{le ext Ber}
Let $\gG$ be a finite group with the identity element $e$.
For any $ b \colon I \raro \gG $ denote
\begin{equation}\label{be(x)}
  B\pn (x) \;:=\; b(x_1) \cdot b(x_2) \cdot \ldots \cdot b(x_n)
  \;\;,\;\;  x ={\{x_n\}}_{n=1}^{\infty} \in X_\rho
\end{equation}
and
\begin{equation}\label{gwb(x)}
 \gw_b(x) \;:=\; min \{ n \in \N \;\colon\; B\pn (x) = e \}
 \;\;,\;\; x \in X_\rho  \;.
\end{equation}
Then the transformation $T_\rho^{\gw_b}$,  defined by
$$
X_\rho \ni x \raro  T_\rho^{\gw_b} x := {T_\rho}^{\gw_b(x)}x \in
X_\rho \;,
$$
is an ergodic endomorphism of $X_\rho$, which is in fact a
one-sided Bernoulli shift.
\end{lemma}
\proof Consider the $\gG$-extension of the graph $(I,\rho)$
generated by $b$.

Namely, set $\; {\ti I}_b :=  \; I \times  \gG \;$ and
$\; {{\ti I}_b}\po0 := \gG \;$
with
$$
 s(\ti i) = (s(i) , \gga) \;,\; t(\ti i) = (t(i) , b(i) \cdot \gga)
 \;,\; p(\ti i) = \rho(i)  \;,\; p\po0 (\gga) = {|\gG|}^{-1}
$$
The skew product endomorphism
$ {\ti T}_{\rho,b} $ corresponding to the stochastic graph
$ {\ti I}_b$, acts on the space $X_\rho \times \gG $ by
$$
  {\ti T}_{\rho,b} (x , \gga) = (T_\rho x \;,\; B(x) \cdot \gga)
  \;\;,\;\; x ={\{ x_n \}}_{n=1}^{\infty} \in X_\rho
  \;,\; \gga \in \gG \;.
$$
where $\; B(x) := {b(x_1)}^{-1} \;$. The skew product $ {\ti
T}_{\rho,b} $ can be identified with the  Markov shift $T_{{\ti
I}_b}$ (see Remark \ref{re GSP shift} ). Under this identification
the partition $\gz^{(0)}_{{\ti I}_b}$ coincides with the partition
$$
 \gz\po0 = \nu_{X_\rho} \times \gep_{\gG} =
 \{\ti E(\gga) \}_{\gga \in \gG } \;,
$$
where
$$
 \ti E(\gga) :=
 X_\rho \times \{\gga\} \subseteq  X_\rho \times \gG
 \;\;.\;\; \gga \in \gG \;.
$$
For any $\gga \in \gG$ consider the endomorphism
$ ({\ti T}_{\rho,b})_{\ti E(\gga)} $
induced by $ {\ti T}_{\rho,b} $ on the set $ \ti E(\gga) $. Let
$$
\gf_{\ti E(\gga)} \;\colon\; \ti E(\gga) \ni (x,\gga) \raro
 \gf_{\ti E(\gga)}(x,\gga) \in \N
$$
be the corresponding return functions (\ref{ret fun}).

Since we use the left shifts on $\gG$ in the definition of the skew
product $ \ti T_b $ and they commute with the right shifts, we have
$$
 \gf_{\ti E(\gga)}(x,\gga) =
 \gf_{\ti E(\gga \cdot \gb)}(x,\gga \cdot \gb)
 \;,\; \gga,\gb \in \gG \;,\; x \in X_\rho \;.
$$
Hence with (\ref{gwb(x)}) and (\ref{be(x)}) we have
$$
  \gw^b(x) = \gf_{\ti E(\gga)}(x,\gga) \;\;,\;\;
$$
and
$$
 {{\ti T}_b}^{\gw^b(x)}(x,\gga) = (T^{\gw^b(x)}x,\gga)
 \;.\; \gga \in \gG \;,\; x \in X_\rho \;,
$$
Thus $T^{\gw^b}$ is isomorphic to the endomorphisms $(T_{{\ti
I}_b})_{D(\gga)}$ induced by the Markov shift $T_{{\ti I}_b}$ on
elements $D(\gga)$ of the partition $ {\gz^{(0)}_{{\ti I}_b}}$. So
that $T^{\gw^b}$ is a Bernoulli shift by Proposition \ref{Tu}.
\qed

\medskip

{\bf Proof of Theorem \ref{th rho cohom} } For given two functions
$a_1$ an $a_2$ put
$$
 b \colon I \ni i \raro b(i) :=
 (a_1(i),a_2(i)) \in \gG := \Ad \times \Ad \;.
$$
and denote for $k=1,2$
$$
 A_k^{\gw_b}(x) :=
 A_k(T_\rho^{\gw_b(x)-1}x) \cdot \ldots \cdot A_k(T_\rho x)
 \cdot A_k(x)   \;\;,\;\; k=1,2
$$
with $A_1$ and $A_2$  defined by (\ref{Am(x)}).
Then by definition of $b$ and $\gw_b$ we have
$$
  A_2^{\gw_b}(x) := A_1^{\gw_b}(x) = e  \;\;,\;\; x \in X_\rho \;,
$$
where $e$ is the identity of $\Ad$. The equality
(\ref{eq rho cohom}) implies
$$
A_2^{\gw_b}(x) \cdot W(x) = W(T_\rho^{\gw_b} x) \cdot A_1^{\gw_b}(x)
$$
and then $\; W(T_\rho^{\gw_b} (x) = W(x) \;$ a.e. on $X_\rho$.
By Lemma \ref{le ext Ber} $\; T_\rho^{\gw_b} $ is ergodic and
hence $W(x)$ is constant a.e. on $X_\rho$. \qed


\bigskip

\subsection{Equivalent extensions.}
\label{ss4.3}

Let $d \in \N$, and $H$ be an irreducible positively recurrent
stochastic graph. Given a function $\; a \colon H \ni h \raro a(h)
\in \Ad \;$ consider the graph skew product $d$-extension $\bH_a$
of $H$ generated by the function $a$ (See Example \ref{GSP}).
Recall that the skew product endomorphism $ \bT_{H,a} $,
corresponding to $ \bH_a $, acts on the space $X_H \times Y_d $ by
$$
  \bT_{H,a} (x,y) = (T_H x , A(x) y)
  \;\;,\;\; x ={\{ x_n \}}_{n=1}^{\infty} \in X_H \;,\; y \in Y_d \;.
$$
where $\; A(x) := {a(x_1)}^{-1} \;$.
We shall use Definition \ref{cohom}
\begin{theorem}\label{Equ ext} Let
$\; \pi_k \colon \bH_{a_k}  \raro  H \;,\; k = 1,2 \;$, be two
$d$-extensions of $H$ generated by functions $a_1$ and $a_2$,
respectively, and let the functions $\; A_k \colon X_H \raro \Ad
\;,\; k=1,2 \;$ are defined by
 \begin{equation}\label{Ak = ak -1}
  A_k(x) := a_k(x_1)^{-1}
  \;\;,\;\; x ={\{ x_n \}}_{n=1}^{\infty} \in X_H \;.
 \end{equation}
 Then the following two conditions are equivalent
 \begin{enumerate}
  \item[(i)] $A_1$ and $A_2$ cohomologous with respect to
   $T_H$, i.e. there exists a
   measurable map  $\; W \colon X_H \raro \Ad \;$ such that
   \begin{equation}\label{A coho}
    A_2(x) \cdot W(x) = W(Tx) \cdot A_1(x) \;\;,\;\; x \in X_H \;,
   \end{equation}
  \item[(ii)] $a_1$ and $a_2$ cohomologous with respect to $H$, i.e.
   there exists a map $\; w \colon H\po0 \raro \Ad \;$ such that
   \begin{equation}\label{a coho}
    a_2(h) \cdot w(s(h)) = w(t(h)) \cdot a_1(h) \;,\; h \in H
   \end{equation}
 \end{enumerate}
\end{theorem}
\proof It is obvious that (\ref{a coho}) implies (\ref{A coho}) with
\begin{equation}\label{W(x)}
 W(x) = w(t(x_1)) \;\;,\;\; x ={\{ x_n \}}_{n=1}^{\infty} \in X_H
\end{equation}
That is (ii) implies (i).

To prove the converse, suppose that (\ref{A coho}) holds with a
suitable measurable function $\; W \colon X_H \raro \Ad \;$.

We have to show that the function $W(x)$ necessarily has the form
(\ref{W(x)}), i.e. $W(x)$ is constant on each element $D(u) =
Z_1^{-1}(_uH)$ of the partition $\; \gz_H\po0 = \{ D(u) \;,\; u
\in H\po0 \} \;$.

To this purpose we shall use induced endomorphisms, which are
defined as follows.

Fix an atom $D(u)$ of the partition $ \gz_H\po0 $ and consider the
endomorphism $ T_u := (T_H)_{D(u)} $, induced by the shift $T_H$
on $D(u)$, see Section \ref{ss3.1}. In accordance with the general
definition (\ref{ret fun}), the return function
$$
 \gf_u(x)=\gf_{D(u)}(x) := min \{ n \geq 1 \colon T_H^{n}x \in D(u) \}
 \;\;,\;\; x \in D(u)
$$
induces $T_u$ by $\; T_ux = T_H^{\gf_u(x)}x \;$. By Proposition
\ref{Tu} the induced endomorphism $T_u$ is isomorphic to the
Bernoulli shift $T_{\rho_u}$, where $\; I_u= \bigcup_{n=1}^\infty
I_{u,n} \;$, and $\; \rho_u =  \{\rho_u(i)\}_{i \in I_u} \;$ are
defined by (\ref{I u n}) and (\ref{rho u}).
That is, $I_{u,n}$ consists of all $\; h_1h_2 \ldots h_n \in H\pn
\;$ such that
$$
 t(h_1) = s(h_n) = u \;\;,\;\; s(h_m) = t(h_{m+1}) \neq u \;,\;
  m = 1,2, \ldots ,n-1
$$
and
$$
\rho_u(i) = p(h_1)  p(h_2) \ldots  p(h_n) \;\;,\;\;
        i = h_1h_2 \ldots h_n \in I_{u,n} \;.\; n \in \N \;.
$$
For any $u \in H\po0$ and $k=1,2$ set
$$
  A_k^{\gf_u}(x) := A_k(T^{\gf_u(x)-1}x)
  \cdot \ldots \cdot A_k(Tx) \cdot A_k(x)
   \;\;.\;\; x \in D(u)
$$
and
$$
 a_k^u (i) := a_k(h_1) \cdot a_k(h_2) \cdot \ldots \cdot a_k(h_n)
 \;\;,\;\;  i = h_1h_2 \ldots h_n \in I_{u,n} \;.
$$
Then
$$
 A_k^{\gf_u}(x) = {a_k^u(i)}^{-1} \;\;,\;\;
 x \in B_u(i) \subseteq D(u) \;.\; i \in I_u \;,
$$
where
$$
  B_u(i) := \{ x ={\{ x_n \}}_{n=1}^{\infty} \in D(u) \;\colon\;
   (x_1x_2 \ldots x_n) = i \in I_{u,n} \;.
$$
Then the equality (\ref{A coho}) implies
\begin{equation}\label{Au cohom}
 A_2^{\gf_u}(x) \cdot W(x) =W(T_u x) \cdot A_1^{\gf_u}(x)
 \;\;.\;\; x \in D(u) \;,
\end{equation}
i.e. $A_1^{\gf_u}$ and $A_2^{\gf_u}$ are cohomologous on $D(u)$
with respect to $ T_u = T_H^{\gf_u} $.

Since for any fix $u \in D(u)$ the partition $\; \gz_u = \{ B_u(i)
\}_{i \in I_u} \;$ is a one-sided Bernoulli generator for $T_u$,
we may apply Theorem \ref{th rho cohom} with the Bernoulli shift $
T_u = T_{\rho_u}$ and with the functions $a_k^u \;,\; k=1,2$.
Therefore, it follows from (\ref{Au cohom}) that there exists $\;
w(u) \in \Ad \;$ such that $\; W(x)=w(u) \;$ a.e. on $ D(u) $.

For every $u$ we have now an element $w(u)$ such that $ W(x) =
w(u) = w(t(x_1))$ for a.e. $ x \in D(u) $.
Hence $ W(x) $ is of the form (\ref{W(x)}), $\; W(T_H x) =
w(s(x_1)) \;$. Thus (\ref{A coho}) implies (\ref{a coho}). \qed

\medskip

As a consequence we obtain

\begin{theorem}\label{Equ ext 1} Let
$\; \pi_k \colon \bH_{a_k}  \raro  H \;,\; k = 1,2 \;$, be two
$d$-extensions of $H$ generated by functions $a_1$ and $a_2$,
respectively.
Let also  $\psi \colon H \raro I $ be an
homomorphism of degree $1$. Suppose $ d = d(\bT_{H,a_1}) =
d(\bT_{H,a_2}) $. Then the endomorphisms $\bT_{H,a_1}$ and
$\bT_{H,a_2}$ are isomorphic iff $a_1$ and $a_2$ cohomologous with
respect to $H$.
\end{theorem}
\proof Since $d(\psi) = 1$ the factor map $\Psi := \Phi_\psi
\colon X_H \raro X_I$ is an isomorphism. Consider two skew
products over the Bernoulli shift $T_\rho$
$$
 \bT_k (x,y) = ( T_\rho x , B_k(x) y)
 \;\;,\;\; (x,y) \in X_\rho \times Y_d
 \;\;,\;\; k = 1,2
$$
where $ B_k(x) := A_k(\Psi^{-1}x) $ and $A_k$ induced by $a_k$ as
above (\ref{Ak = ak -1}). Each of the shifts $T_{\bH_{a_k}}$ is a
simple $\rho$-uniform endomorphism by Theorem \ref{simp mar}. The
skew products $\bT_{H,a_k}$ as well as the shifts $T_{\bH_{a_k}}$,
are isomorphic to $\bT_k$.
They are $\rho$-uniform endomorphisms
and $d = d(\bT_1) = d(\bT_2).$ By Theorem \ref{simple} $\bT_1$ and
$\bT_2$ are isomorphic iff the functions $B_1$ and $B_2$ are
cohomologous with respect to $T_\rho$. This means that the
functions $A_1$ and $A_2$ are cohomologous with respect to $T_H$.
Finally, by Theorem \ref{Equ ext} the last condition holds iff
$a_1$ and $a_2$ cohomologous with respect to $H$. \qed


\bigskip

\subsection{ $\bldGSP$-extensions and persistent $\bldd$-partitions.}
\label{ss4.4}

Let $H$ be a stochastic graph and $(I,\rho)$ a standard Bernoulli
graph. In this section we study  extensions of the form
\begin{equation}\label{pi psi}
(\pi,\psi) \;\colon\; \xymatrix{ \bH \ar[r]^{\pi} & H \ar[r]^{\psi} & I}
\end{equation}
where the graph $H$ be an extension of the Bernoulli graph
$(I,\rho)$ by a homomorphism $\psi$ of degree $ d(\psi) = 1 $ and
$ \bH = \bH_a $ be a graph skew product $d$-extension of $H$
generated by a function $ a \colon H \raro \Ad $ (See Example
\ref{GSP}). The diagrams of the above form (\ref{pi psi}) will be
referred to $\bldl \bldpi , \bldpsi \bldr${\bf -extensions}. We
assume that the graph $ \bH $ is irreducible, i.e. the
corresponding Markov shift $T_\bH$ and skew product $\bT_{H,a} $
are ergodic.

\medskip

Fixing an extension (\ref{pi psi}) and setting $\; J = H\po0 \;$,
we identify $H$ with $ I \times J $ such that
$$
 \psi (h) = i \;,\; s(h) = j \;,\; t(h) = f_i (j)
 \;,\; p(h) = \rho(i)
$$
for any $ h = (i,j) \in H = I \times J $.
Here the maps $ f_i \colon J \raro J $ are uniquely determined by
$$
 f_ij = t(i,j) \;\;,\;\; (i,j) \in  I \times J
$$
and the semigroup $ \cS (\psi) $, generated by $ \{ f_i \;,\; i
\in I \} $ is $1$-contractive, since $ d(\psi) = 1 \;$ (Theorem
\ref{d(phi)}).

The $d$-extension $ \bH = \bH_a $ is described now as follows:
\begin{equation}\label{bar H}
 \bH = I \times J \times Y_d \;\; , \;\;
 \bH\po0 = H\po0 \times Y_d = J \times Y_d \;,
\end{equation}
where for any $ \bh = (i,j,y) $
\begin{equation}\label{stp bar H}
  s(\bh) = (j,y)   \;,\; t(\bh) = (f_ij , a(i,j)y) \;,\;
  p(\bh) = \rho(i) \;,\; a(h) = a(i,j)
\end{equation}
The homomorphisms $\psi$, $\pi$ and $\phi = \psi \circ \pi $ are
defined by
$$
 \pi (\bh) = h = (i,j) \;\;,\;\; \pi\po0 (j,y) = j \;\;,\;\;
 \phi (\bh) = \psi (h) = i \;,
$$
where $ d(\phi) = d(\pi) = d $ and the diagram
\begin{equation}\label{diag I J Yd}
\xymatrix{
 I \times J \times Y_d \ar[d]^{\pi}  \ar[dr]^{\phi} &   \\
 I \times J \ar[r]^{\psi}                           & I  }
\end{equation}
commutes.
The semigroup $\cS(\phi)$ can be described now as a $d$-extension
$\bcS = \bcS (\pi,\psi)$ of the semigroup $\cS(\psi)$.

Set
\begin{equation}\label{bfi (j,y)}
 \baf_i(j,y) := t(i,j,y) = (f_i j \;,\; a(i,j) y)
 \;,\; (j,y) \in J \times Y_d \;,\; i \in I  \;.
\end{equation}
The maps $ \baf_i $ act on $J \times Y_d $.

The semigroup $\bcS$, generated by $\; \baf_i \;,\; i \in I \;$,
consists of all maps of the form:
$$
 \{ \baf_{i_1i_2 \ldots i_n} \;,\; i_1i_2 \ldots i_n \in I^n \;,\;
   n \in \N \} \;,
$$
where
$$
 \baf_{i_1i_2 \ldots i_n}(j,y) = (f_{i_1i_2 \ldots i_n} j ,
 a(i_1i_2 \ldots i_n,j) y)
$$
and
$$
 a(i_1i_2 \ldots i_n,j) := a(i_1,f_{i_2i_3 \ldots i_n} j )
 \ldots a(i_{n-1},f_{i_n} j ) a(i_n,j) \;.
$$
Note that $\; \cS(\psi) \ni f \raro \baf \in \bcS \;$ is an
isomorphism between the semigroups.
\begin{proposition}\label{S(phi) bar H}
The semigroup $\bcS $ is $d$-contractive and
its persistent $d$-sets are of the form
$$
 \cL (\bcS) = \{ L_j , j \in J \} \;\;,\;\; L_j := \{ j \} \times Y_d  \;.
$$
\end{proposition}
\proof By Theorem \ref{d(phi)} the semigroup $ \bcS $ is
$d$-contractive and $\cS (\psi)$ is $1$-contractive, since
$ d(\phi) = d $ and $ d(\psi) = 1 $.

For any finite set $ F \subseteq J \times Y_d $ the set $ E :=
\pi\po0 (F) \subseteq J $ is also finite. Since the semigroup
$\cS(\psi)$ is $1$-contractive there exist $ i_1i_2 \ldots i_n \in
I^n $ and $j \in J$ such that $\; f_{i_1i_2 \ldots i_n}(E) = \{j\}
\;$ and hence $\; \baf_{i_1i_2 \ldots i_n}(F) \subseteq L_j \;$.
On the other hand $ d = |L_j| = |\baf_i(L_j)| $ for all $i \in I ,
j \in J$. Thus the sets $L_j$ and only they are persistent sets
for the semigroup $\bcS$. \qed

\medskip

For every $E \subseteq J$ set $ \bE := {\pi\po0}^{-1}E = E \times
Y_d $.
\begin{definition}\label{pers part}
Let the semigroup $\bcS$ be as above.
 \begin{enumerate}
  \item[(i)] A subset $ R $ of $ \bH\po0 = J \times Y_d $ will be
   called {\bf transversal} with respect to $\pi\po0$
   if $ \pi\po0 (R) = H\po0 = J $ and the restriction
   $ \pi\po0 {|}_R \colon R \raro J $ is a bijection.
   A partition $\; r = \{ R_1,R_2, \ldots ,R_d \} \;$
   will be called {\bf transversal} with respect to
   $\pi\po0$ if all the set $R_1,R_2, \ldots ,R_d$ are transversal.
  \item[(ii)] A transversal partition $r$ will be called {\bf
   persistent} with respect to semigroup $ \bcS $, if for every
   transversal partition $r_1$ and
   every finite subset $E \subseteq J$ there exists
   $ \baf \in \bcS $ such that
   $\; \baf^{-1}r_1 \;{|}_\bE  = r \;{|}_\bE \;.$
 \end{enumerate}
\end{definition}
Denote by $\cR$ the set of all transversal partitions and by $\cR
(\bcS)$ the set of all persistent partitions for the semigroup
$\bcS $.

If a set $R$ is transversal then for any $i \in I$ the set
$\baf^{-1}R$ is transversal. Hence for any $r = \{ R_1,R_2, \ldots
,R_d \} \in \cR $ and $ \baf \in \bcS $ we have
$$
 \baf^{-1} r := \{ \; \baf^{-1}R_1 , \baf^{-1}R_2, \ldots ,
 \baf^{-1}R_d \; \} \in \cR \;.
$$
Further we shall use this action
$$
\cR \ni r \; \raro \;  \baf^{-1} r \in \cR \;\;,\;\; \baf \in \bcS
$$
of the semigroup $\bcS$ on $\cR$. The following lemma shows that
$ \cR(\bcS)$ is an attracting set for $\cR$ with respect to the
action in a natural sense.
\begin{lemma}\label{gL(phi)}
 \begin{enumerate}
  \item[(i)] The set $\cR(\bcS)$ is not empty.
  \item[(ii)] $\baf^{-1} \cR (\bcS) \subseteq  \cR (\bcS)$ for all
  $\baf \in \bcS$
  \item[(iii)] $\cR (\bcS)$ is the least subset of $\cR$ with the
   property (ii).
 \end{enumerate}
\end{lemma}
\proof  Consider a subset $\cR_0 (\bcS)$ of $\cR$ consisting of all
$r \in \cR$ having the following property:
 \begin{enumerate}
  \item[]
   For any finite subset $E \subseteq J $ there exists
   $\; f \in \cS(\psi) $ such that $f(E)$ is a singleton, i.e.
   $|f(E)|=1$, and $ r \;{|}_\bE =  \baf^{-1}\gep \;{|}_\bE \;$,
   where $ \bE := E \times Y_d $ and $ \gep = \gep_{J \times Y_d} $.
 \end{enumerate}
We show that $\cR (\bcS) = \cR_0(\bcS) \neq \emptyset$.

Take a sequence $\; E_n \nearrow J \;,\; |E_n| < \infty \;$.
Since $d(\cS(\psi)) = 1 \;$ we can find a sequence $g_n \in \cS(\psi)$ such that
for all $n \in \N$ and $ f_n := g_n \cdot \ldots \cdot g_2 \cdot g_1$, the set
$f_n (E_n)$ is single-point, i.e. $|f_n (E_n)|=1$.
Using the decreasing sequence of partitions
$$
 \gep \geq \baf_1^{-1} \gep \geq \baf_2^{-1} \gep \geq \; \ldots \;
 \geq \baf_n^{-1} \gep \geq \ldots
$$
set $\; r_0 := \bigwedge_{n=1}^{\infty} {\baf_n}^{-1}\gep \;$.
Since $|f_n (E_n)|=1$ the restriction $r_0 \;{|}_{\bE_n} $
consists of $d$ sets, whose projections on $J$ are $E_n$. Hence $
r_0 \in \cR$ and $r_0 \;{|}_{\bE_n} = \baf_n^{-1}\gep
\;{|}_{\bE_n} $. We see that $r_0 \in \cR_0(\bcS)$, i.e. $\cR_0$
is not empty,

Let $r \in \cR_0(\bcS)$ and $r_1 \in \cR$. For any finite subset
$E \subseteq J$ there exists $\baf \in \cR(\bcS)$ such that $ r
\;{|}_\bE =  \baf^{-1}\gep \;{|}_\bE \;$. Then
$$
 \baf^{-1}r_1 \;{|}_\bE \leq \baf^{-1}\gep \;{|}_\bE = r \;{|}_\bE \;,
$$
Since each of the partitions consists of $d$ elements, we have also
$\baf^{-1}r_1 \;{|}_\bE = r \;{|}_\bE$.
Hence $r \in \cR(\bcS)$.

Conversely, let $r \in \cR(\bcS)$ and $E$ be a finite subset of
$J$. There are exist $\baf \in \bcS$ and $r_1 \in \cR$ such that $
r_1 \;{|}_\bE = \baf^{-1}\gep \;{|}_\bE \;$.
On the other hand, since $r \in \cR(\bcS)$, we can choose $\baf_1
\in \bcS$ for which $ \baf_1^{-1}r_1 \;{|}_\bE = r \;{|}_\bE $.
Hence
$$
 r \;{|}_\bE = \baf_1^{-1}r_1 = \baf_1^{-1}\baf^{-1}\gep \;{|}_\bE
 = (\baf^{-1}\baf_1)^{-1}\gep \;{|}_\bE.
$$
We see that $r \in \cR_0(\bcS)$ and thus $\cR (\bcS) = \cR_0$ and
Part $(i)$ follows.

Parts $(ii)$ and $(iii)$ follow in the same manner by the
definition of $\cR (\bcS)$ and by the equality $\cR (\bcS) =
\cR_0(\bcS)$. \qed


\bigskip

\subsection{ Irreducible $\bldd$-extensions.}
\label{ss4.5}

In this section we continue to study $(\pi,\psi)$-extensions of
the form (\ref{pi psi})
$$
 (\pi,\psi) \;\colon\; \xymatrix{
 \bH \ar[r]^{\pi} & H \ar[r]^{\psi} & I } \;,
$$
where  the graph $H$ is an extension of the Bernoulli graph
$(I,\rho)$ by a homomorphism $\psi$ of degree $ 1 $ and $ \bH =
\bH_a $ be a GSP $d$-extension of $H$, generated by a function $ a
\colon H \raro \Ad $.

Fix $d$ and $(I,\rho)$ and consider the set $\; \Ext \;$ of all
$(\pi,\psi)$-extensions of the form (\ref{pi psi}). This set is
equipped with a natural partial order and with an equivalence
relation as follows
\begin{definition}\label{partial order} Let
$\; (\pi,\psi) \colon
\xymatrix{ \bH \ar[r]^{\pi} & H \ar[r]^{\psi} & I } \;$
and
$\; (\pi_1,\psi_1) \colon
      \xymatrix{\bH_1 \ar[r]^{\pi_1} & H_1 \ar[r]^{\psi_1} & I } \;$
be two $(\pi,\psi)$-extensions from $\Ext$. Let $a$ and $a_1$ be
the functions, which generate the extensions $\bH$ and $\bH_1$,
respectively.
\begin{enumerate}
 \item[(i)] A homomorphism $\; \bgk \colon \bH \raro \bH_1 \;$
  is said to be a {\bf trivializable} $d$-extension of a homomorphism
  $\; \gk \colon H \raro H_1 \;$, if the square part of
  Diagram \ref{diag order} ( below ) commutes and the functions

  $a_1 \circ \chi$ and $a$ are cohomologous with respect to $H$.
 \item[(ii)] We shall say that
  $$
      (\pi_1,\psi_1) \; \preceq \; (\pi,\psi)
  $$
  if there is a commutative diagram
  \begin{equation}\label{diag order}
   \xymatrix{
    \bH \ar[d]^{\bgk} \ar[r]^{\pi}
  & H \ar[d]^{\gk} \ar[dr]^{\psi} & \\
    \bH_1 \ar[r]^{\pi_1}
  & H_1 \ar[r]^{\psi_1}           & I }\;,
  \end{equation}
  where $\; \bgk \in \Hom (\bH,\bH_1) \;$
  is a trivializable $d$-extension of $\; \gk \in \Hom (H,H_1) \;$.
 \item[(iii)]  We shall say that
  $$
      (\pi_1,\psi_1) \; \sim \; (\pi,\psi) \;
  $$
  if there is commutative Diagram \ref{diag order},
  where both $\; \gk \;\colon\; H \raro H_1 \;$ and its $d$-extension
  $\; \bgk \;\colon\; \bH \raro \bH_1 \;$ are isomorphisms.
 \end{enumerate}
\end{definition}
In connection with Part (i) of the definition, note that an
extension $\; \bgk \;\colon\; \bH \raro \bH_1 \;$ is trivializable
iff it is equivalent to a trivial extension of $\; \gk \;:\; H
\raro H_1 \;$ (see Remark \ref{triv exten}).

It can be checked also that $\; (\pi_1,\psi_1) \preceq (\pi,\psi)
\;$ and $\; (\pi,\psi) \preceq (\pi_1,\psi_1) \;$ imply $\;
(\pi_1,\psi_1) \sim (\pi,\psi) \;$, but we do not use the fact in
this paper.

\medskip

Our aim now is to describe "minimal" elements of $\; (\; \Ext
\;,\; \preceq \;) \;$.
\begin{definition}\label{irreduc}
An extension $(\pi,\psi) \in \Ext$ is called {\bf irreducible} if
$(\pi_1,\psi_1) \sim (\pi,\psi)$ as soon as $ (\pi_1,\psi_1) \in
\Ext$ and $(\pi_1,\psi_1) \preceq (\pi,\psi)$.
\end{definition}
\begin{theorem}\label{reduc ext}
For any $\; (\pi,\psi) \in \Ext \;$ there exists a unique up to
equivalence irreducible $(\pi,\psi)$-extension
$\; (\pi_*,\psi_*) \in \Ext \;$ such that
$\; (\pi_*,\psi_*) \preceq (\pi,\psi) \;$.
\end{theorem}
To prove the theorem we fix a pair $(\pi,\psi) \in \Ext$ and again
use the identification (\ref{bar H}). Namely,
\begin{equation}\label{diagr I J Yd}
 (\pi, \psi) \;\;:\;\;
 \xymatrix{
   \bH = I \times J \times Y_d \ar[r]^-\pi
 & H = I \times J \ar[r]^-\psi & I }
\end{equation}
where $H\po0 = J$ and $\bH\po0 = H\po0 \times Y_d = J \times Y_d $
as in Section \ref{ss4.4}.

We construct the desired irreducible $(\pi_*,\psi_*)$-extension
and a corresponding commutative diagram
\begin{equation}\label{diag reduc ext}
 \xymatrix{
    \bH \ar[d]^{\bgk_*} \ar[r]^{\pi}
  & H \ar[d]^{\gk_*} \ar[dr]^{\psi} &   \\
    \bH_* \ar[r]^{\pi_*}
  & H_* \ar[r]^{\psi_*}             & I  }
\end{equation}
by means of the semigroup $ \bcS = \bcS (\pi,\psi) $ and its
persistent partitions $\cR (\bcS)$.
\begin{definition}\label{reduc part}
A partition $\xi$ of $J=H\po0$ is called {\bf reducing} partition if
the following two conditions hold
 \begin{enumerate}
  \item[(i)] $f^{-1}\xi \leq \xi$ for all $f \in \cS (\psi)$, i.e.
  $\xi$ is $\cS(\psi)$-invariant
  \item[(ii)] For any element $C \in \xi$ denote $\bC := \pi^{-1}C$
   and let  $r {|}_\bC $ be the restriction of the partition
   $r \in \cR (\bcS)$ on the set $\bC$.
   Then  all the  partitions $\; r {|}_\bC \;,\; r \in \cR (\bcS) \;$
   coincide with each other.
 \end{enumerate}
\end{definition}
Consider the set $\Xi$ of all reducing partitions $\xi$ on $H\po0$.

For any $\xi \in \Xi$ we have
$\; {\pi\po0}^{-1}\xi = \xi \times \nu_{Y_d} \;$ and the partition
$\; {\pi\po0}^{-1}\xi \vee r \;$ does not depend on the choice of
$r \in \cR(\bcS)$.
So that we may set
\begin{equation}\label{bxi}
  \bxi := {\pi\po0}^{-1}\xi \vee r \;\;,\;\; \xi \in \Xi
\end{equation}
and $\bXi := \{\bxi \;\colon\; \xi \in \Xi \}$ on $\bH\po0$.

Since $\xi$ is $\cS(\psi)$-invariant and $\cR(\bcS)$ is
$\bcS$-invariant by Lemma \ref{gL(phi)}, the partition $\bxi$ is
also $\bcS$-invariant.

Therefore we may introduce the {\bf factor pair}
\begin{equation}\label{factor pair}
 \xymatrix{
   \bH{/}_\bxi \ar[r]^{\pi{/}_\xi}
 & H{/}_\xi \ar[r]^{\psi{/}_\xi} & I }
\end{equation}
Namely, we set
$$
 H{/}_\xi := I \times J{/}_\xi \;\;,\;\;
 \bH{/}_\bxi := I \times J{/}_\xi \times Y_d
$$
Any element of $\bxi$ consists of $d$  elements of the form $\;
R_y^C \;,\; y \in Y_d \;$, where $C \in \xi$ and $\pi\po0(R_y^C) =
C$. Hence. by possibly passing to an equivalent extension, we may
assume that $R_y^C = C \times \{y\}$, i.e. $\bxi = \xi \times
\gep_{Y_d}$. This means that the function $a = a(i,j)$, generating
the extension $\bH = \bH_a$, does not depend on $j$ on the
elements of $\xi$.
Hence the equalities (\ref{stp bar H}) and
(\ref{bfi (j,y)}) well define $a{/}_\xi$ and $ \bH{/}_\bxi :=
(H{/}_\xi)_{a{/}_\xi}$. Thus we have shown
\begin{proposition}\label{natur proj}
For any $\xi \in \Xi$ the natural projections
$$
 {\pi_\bxi}\po0 \colon \bH\po0 \raro \bH\po0{/}_\bxi \;\;,\;\;
 {\pi_\xi}\po0 \colon H\po0 \raro H\po0{/}_\xi
$$
uniquely determine $(\pi{/}_\xi,\psi{/}_\xi) \in \Ext$ such that
$ (\pi{/}_\xi,\psi{/}_\xi) \preceq (\pi,\psi) $ with the
coresponding commutative diagram
 \begin{equation}\label{diag /xi}
  \xymatrix{
     \bH \ar[d]^{\pi_\bxi} \ar[r]^{\pi}
  & H \ar[d]^{\pi_\xi} \ar[dr]^{\psi} &  \\
     \bH{/}_\bxi \ar[r]^{\pi{/}_\xi}
  & H{/}_\xi \ar[r]^{\psi{/}_\xi}     & I }
 \end{equation}
\end{proposition}
Conversely
\begin{proposition}\label{convers}
For any $\; (\pi_1,\psi_1) \; \preceq \; (\pi,\psi) \;$
there exists $\xi \in \Xi$ such that
$(\pi{/}_\xi,\psi{/}_\xi) \sim (\pi_1,\psi_1)$
\end{proposition}
\proof Take the map $ \gk\po0 \colon H\po0 \raro {H_1}\po0 $
induced by homomorphism $ \gk \colon H \raro H_1 $ from Diagram
\ref{diag ext} and set $\xi := {\gk\po0}^{-1}\gep_{{H_1}\po0}$.
Then $\xi \in \Xi$ and it is desired \qed

\medskip

{\bf Proof of Theorem \ref{reduc ext}.} It is easily to see that
$\Xi$ is a lattice, i.e. $ \xi_1 \vee \xi_2 \in \Xi $ and $ \xi_1
\wedge \xi_2 \in \Xi $ for all $\xi_1 , \xi_2 \;\in\;\Xi $.
Herewith, $\Xi$ has the least element. Denote the least element by
$\xi_* $ and let $ \bxi_* := \overline{(\xi_*)}$ be the
corresponding partition of $\bH\po0$. Note that $ \bxi_*$ is the
least element of $\bXi $.
Herewith $ \bxi_*$ is the least partition of $\bH\po0$ such that
for all $r \in \cR (\bcS)$ and every $C \in \xi$ the restriction
$r \;{|}_\bC $ consists precisely of $d$ elements.

Putting $\xi = \xi_*$ in Diagram \ref{diag /xi} (Proposition
\ref{natur proj}) we obtain Diagram \ref{diag reduc ext} with
$$
 H_* =   H{/}_{\xi_*}    \;\;,\;\; \bH_* = \bH{/}_{\bxi_*}
 \;\;,\;\; \pi_* = \pi{/}_{\xi_*}   \;\;,\;\; \psi_* =
 \psi{/}_{\xi_*} \;.
$$
and $ (\pi_*,\psi_*) \preceq (\pi,\psi) $.

Using by the above propositions and Lemma \ref{gL(phi)}, we see
that the pair $(\pi_*,\psi_*)$ is irreducible and that it is the
only (up to equivalence) irreducible pair majorized by $
(\pi,\psi) $. \qed
\begin{remark}\label{irred ext}
The above arguments show that a pair $ (\pi,\psi) \in \Ext $ is
irreducible iff $ (\pi_*,\psi_*) = (\pi,\psi) $, i.e. iff $ \xi_*
= \gep_{H\po0}$. The last equality means that the persistent
partitions $ \cR (\bcS) $ separate the points of $H\po0$ in the
following sense: for every pair $ u_1,u_2 \in H\po0 $ there exist
$R_1 \in r_1 \in \cR (\bcS)$ and $R_2 \in r_2 \in \cR (\bcS)$ such
that
$$
 {\pi\po0}^{-1}(u_1) \cap R_1 \cap R_2 \neq \emptyset \;\;,\;\;
 {\pi\po0}^{-1}(u_2) \cap R_1 \cap R_2 = \emptyset \;.
$$
\end{remark}


\bigskip

\section{Canonical form and classification.}
\label{s5}

\bigskip

\subsection{ Main Theorems.}
\label{ss5.1}

The following two theorems claim the existence and uniqueness of
the canonical form of $\rho$-uniform one-sided Markov shifts.
\begin{theorem}\label{canon form}
Let $G$ be a $\rho$-uniform stochastic graph, which is irreducible and
positively recurrent.
Then there exists a $(\pi,\psi)$-extension
\begin{equation}\label{pi,psi}
 (\pi,\psi) \;\colon\;
 \xymatrix{ \bH \ar[r]^{\pi} & H \ar[r]^{\psi} & I }
\end{equation}
such that
 \begin{enumerate}
  \item[(i)] The shifts $T_G$ and $T_\bH$ are isomorphic,
  \item[(ii)] $(\pi,\psi) \in \Ext$, where $d = d(T_G)$ is the minimal
  index of the shift $T_G$,
  \item[(iii)] The extension $(\pi,\psi) $ is irreducible.
 \end{enumerate}
\end{theorem}
\proof Combining the results of Theorems \ref{reduc ext} and
\ref{phi bar d(T)} we obtain from Diagrams \ref{diag reduc ext}
and \ref{diag bar H Gn} the following commuting diagram
\begin{equation}\label{diag main}
   \xymatrix{
          \bH \ar[dd]^{\bgk} \ar[rd]^{\pi} \ar[rr]^{\bpsi}
     &&   G\pn \ar[dd]^{\phi} \ar[rr]^{\pi\pn}
     &&   G  \\
     &    H \ar[d]^{\gk} \ar[dr]^{\psi}       &    &&     \\
          \bH_* \ar[r]^{\pi_*}
     &    H_* \ar[r]^{\psi_*}        &  I &&      }
\end{equation}
Here, $\pi$, $\pi_*$ and $\phi$ are homomorphisms of degree
$d = d(T_G)$, all other homomorphisms are of degree $1$, and the
extension $(\pi_*,\psi_*) \in \Ext $ is irreducible.

Since $G$ and $ \bH_* $ have a common extension $\bH$ of degree
$1$, the shifts $T_G$ and $T_{\bH_*}$ are isomorphic. Thus the the
extension $(\pi_*,\psi_*) $ is desired.
\qed
\begin{definition}\label{def canon form}
We shall say that $T_\bH$ is a {\bf canonical form} of the shift
$T_G$, if there exists an extension (\ref{pi,psi}) satisfying the
conditions of Theorem \ref{canon form}. Herewith the graph $\bH$
is said to be the {\bf canonical graph} for $T_G$.
\end{definition}
Theorem \ref{canon form} states the existence of the canonical
form. Turn to the uniqueness.
\begin{theorem}\label{classification}
Let $G_1$ and $G_2$ be two $\rho$-uniform stochastic graphs, which
are irreducible and satisfy the positive recurrence condition.
Suppose the shifts $T_{G_1}$ and $T_{G_1}$ are represented in the
canonical form $T_{\bH_1}$ and $T_{\bH_2}$, respectively, and let
\begin{equation}\label{pi12,psi12}
  (\pi_k,\psi_k) \;\colon\;
     \xymatrix{
     \bH_k \ar[r]^{\pi_k} & H \ar[r]^{\psi_k} & I }  \;\;,\;\; k=1,2
\end{equation}
be corresponding canonical $(\pi,\psi) $-extensions.

Then the following conditions are equivalent
 \begin{enumerate}
  \item[(i)] The shifts $T_{G_1}$ and $T_{G_2}$ are isomorphic,
   $(T_{G_1} \sim T_{G_2})$.
  \item[(ii)] The graphs $\bH_1$ and $\bH_2$ are isomorphic,
   $(\bH_1 \sim \bH_2)$.
  \item[(iii)] The extensions $(\pi_1,\psi_1) $ and $(\pi_1,\psi_1) $

  are equivalent,$ ( (\pi_1,\psi_1) \sim (\pi_2,\psi_2) ) $.
 \end{enumerate}
\end{theorem}
\proof By the definition we have $T_{G_1} \sim T_{\bH_1}$ ,
$T_{G_2} \sim T_{\bH_2}$ and
$$
 (\pi_1,\psi_1) \sim (\pi_2,\psi_2) \;\; \RRaro \;\; \bH_1 \sim \bH_2
 \;\; \RRaro \;\; T_{\bH_1} \sim T_{\bH_2}
$$
Thus we need to prove only
\begin{equation}\label{RRaro}
 T_{\bH_1} \sim T_{\bH_2} \;\; \RRaro
 \;\; (\pi_1,\psi_1) \sim (\pi_2,\psi_2)
\end{equation}
Suppose $T_{\bH_1} \sim T_{\bH_2} $ and let $\; a_k \colon H_k
\raro \Ad \;,\; k=1,2 \;$, be the functions generating $\bH_k$,
where $\; d = d(T_{\bH_1}) = d(T_{\bH_2}) \;$.

Since both of $ \psi_1 \colon H_1 \raro I $ and $ \psi_2 : H_2
\raro I $ are of degree $1$, we can apply Theorem \ref{degree 1}
and to construct a common extension $H$ of $H_1$ and $H_2$.
Herewith, the corresponding Diagram \ref{diag degree 1} commutes
and the homomorphisms $\; \psi \colon H \raro I \;$ , $\; \chi_1
\colon H \raro H_1 \;$ and $\; \chi_2 \colon H \raro H_2 \;$ are
of degree $1$.

By Remark \ref{triv exten} each of homomorphisms $\; \chi_k \colon
H_{b_k} \raro H_k \;,\; k = 1,2 \;$ admits the trivial extension
$\; \bchi_k \colon \bH_{b_k} \raro \bH_k \;$ with the commuting
diagram
\begin{equation}\label{diag ext k=1,2}
   \xymatrix{
      \bH_{b_k} \ar[d]^{\pi_{b_k}} \ar[r]^{\bchi_k}
    & \bH_k \ar[d]^{\pi_k} \\
        H \ar[r]^{\chi_k}           &  H_k }
\end{equation}
Here $\bchi_k$ is of degree $1$ and $b_k := a_k \circ \chi_k $ for
$k = 1,2$. Since $ d(\bchi_1) = d(\bchi_2) = 1 $ we have
$T_{\bH_1} \sim T_{\bH_{b_1}} $ and $T_{\bH_2} \sim T_{\bH_{b_2}}
$. Therefore $T_{\bH_1} \sim T_{\bH_2} $ implies that the skew
products $\bT_{H,b_1}$ and $\bT_{H,b_2}$ are isomorphic.

Thus we have two GSP $d$-extensions
$\; \pi_{b_k} \colon \bH_{b_k} \raro H \;,\; k = 1,2 \;,\;$
of $H$ and a homomorphism $\psi : H \raro I $ of degree $1$.
Herewith, the number $d$ is the minimal index of $\bT_{H,b_1}$ and
$\bT_{H,b_2}$.
By Theorem \ref{Equ ext 1} the functions $b_1$ and $b_2$ are cohomologous
with respect to $H$.
Hence two constructed $(\pi,\psi)$-extensions
\begin{equation}
  (\pi_{b_k},\psi) \;\colon\;
    \xymatrix{ \bH_{b_k} \ar[r]^{\pi}
       & H \ar[r]^{\psi} & I } \;\;,\;\; k=1,2
\end{equation}
are equivalent, $\; (\pi_{b_1},\psi) \sim (\pi_{b_2},\psi) \;$.

On the other hand by constructing both two diagrams
\begin{equation}
   \xymatrix{
   \bH_{b_k} \ar[d]^{\bgk_k} \ar[r]^{\pi_{b_k}}
 & H \ar[d]^{\gk} \ar[dr]^{\psi} & \\
   \bH_k \ar[r]^{\pi_k}
 & H_k \ar[r]^{\psi_k}    & I } \;\;,\;\; k=1,2
\end{equation}
commute. This means that $\; (\pi_1,\psi_1) \; \preceq \;
(\pi_{b_1},\psi) \;$ and $\; (\pi_2,\psi_2) \; \preceq \;
(\pi_{b_2},\psi) \;$.

The pairs $(\pi_1,\psi_1)$ and $(\pi_2,\psi_2)$ are irreducible
and they are majorized by equivalent pairs. Hence they are
equivalent.

We have shown (\ref{RRaro}).
\qed

\medskip

As a consequence we have also
\begin{theorem}\label{common exten 1}
Under conditions of Theorem \ref{classification} the shifts
$T_{G_1}$ and $T_{G_2}$ are isomorphic iff the graphs $G_1$ and
$G_2$ have a common extension of degree $1$, i.e. there exists a
diagram
\begin{equation}\label{diag com ext 1}
   \xymatrix{
    G_!  &  G \ar[l]_{\phi_1} \ar[r]^{\phi_2}  &  G_2  }
\end{equation}
where homomorphisms $\phi_1$ and $\phi_2$ are of degree $1$.
\end{theorem}
\proof By Theorem \ref{phi bar d(T)} we have two diagram of
 homomorphisms
\begin{equation}\label{diag com ext 2}
 \xymatrix{
 G_k & G_k\pn \ar[l]_{\pi\pn}
 &  \bH_k \ar[l]_{\bpsi_k} \ar[r]^{\pi_k}
 & H_k \ar[r]^{\psi_k} & I   }  \;\;\;\; k=1,2
\end{equation}
where $d(\pi\pn) = d(\bpsi_k) = d(\psi_k) = 1 $ and $\pi_k$ is a
$d$-extension. So that $(\pi_k,\psi_k) \in \Ext$.

By Theorem \ref{reduc ext} each pair $(\pi_k,\psi_k) \;,\; k=1,2
\;$ majorizes an irreducible pair from $ \Ext $. If the the shifts
$T_{G_1}$ and $T_{G_2}$ are isomorphic the irreducible pairs are
equivalent (Theorem \ref{classification}) and we may assume
without loss of generality that they coincide with each other.

Thus there exists $(\pi_0,\psi_0) \in \Ext$ with two commuting
diagrams
\begin{equation}\label{diag reduc ext k=1,2}
 \xymatrix{
    \bH_k \ar[d]^{\bgk_k} \ar[r]^{\pi_k}
  & H_k \ar[d]^{\gk_k} \ar[dr]^{\psi_k} &  \\
    \bH_0 \ar[r]^{\pi_0}
  & H_0 \ar[r]^{\psi_0}    & I  } \;\;\;\; k=1,2
\end{equation}
Passing possibly to equivalent extensions we may also assume that
$\bgk_1$ and $\bgk_2$ are trivial extensions of $\gk_1$ and
$\gk_2$.

By Theorem \ref{degree 1 sharp} and Remark \ref{triv exten} we
find a common extension of degree 1
\begin{equation}\label{diag com ext 3}
   \xymatrix{  H_1  &  H \ar[l]_{\chi_1} \ar[r]^{\chi_2}  &  H_2  }
\end{equation}
of $H_1$ and $H_2$ with the trivial extensions
\begin{equation}\label{diag com ext 4}
   \xymatrix{
   \bH_!  &  \bH \ar[l]_{\bchi_1} \ar[r]^{\bchi_2}  &  \bH_2  }
\end{equation}
of $\chi_1$ and $\chi_2$ such that the corresponding diagram
\begin{equation}\label{diag com ext 5}
   \xymatrix{
 & \bH_2 \ar[dd]_(.75){\pi_2} \ar[rr]^{\bgk_2}&
 & \bH_0 \ar[dd]^{\pi_0} \ar[drrrr]^{\bpsi_0}   & &&& \\
   \bH \ar[dd]_{\pi} \ar[ur]^{\bchi_2} \ar[rr]^(.75){\bchi_1}
     \ar[rrru]^{\bchi}
 && \bH_1 \ar[ur]_{\bgk_1} \ar[dd]^(.25){\pi_1} &&&&& I \\
 &  H_2 \ar[rr]^(.75){\gk_2}&
 &  H_0  \ar[rrrru]^{\psi_0} &&&& \\
    H \ar[ur]^{\chi_2} \ar[rr]_{\chi_1} \ar[rrru]^{\chi}
         & & H_1 \ar[ur]_{\gk_1} & && && \\  }
\end{equation}
commutes. Therefore we have
\begin{equation}\label{diag com ext 6}
 \xymatrix{
  G_1 & G_1\pn \ar[l]_{\pi\pn}  &  \bH_1 \ar[l]_{\bpsi_1}
  & \bH \ar[l]_{\bgk_1} \ar[r]^{\bgk_2} & \bH_2 \ar[r]^{\psi_2}
  & G_2\pn \ar[r]^{\pi\pn} & G_2  }
\end{equation}
Putting $\; G := \bH \;$ and $\; \phi_k := \bgk_k \circ \bpsi_k
\circ \pi\pn \;$ for $\; k = 1,2 \;$, we obtain the desired common
extension of degree $1$ (\ref{diag com ext 1}). \qed


\bigskip

\subsection{Consequences and examples.}
\label{ss5.2}

Consider some particular cases.

\bigskip

\noindent {\bf Extensions of Bernoulli graphs.} Let $(I,\rho)$ be
a standard Bernoulli graph and let $ d \in \N $. Let  $\; a \colon
I \raro \Ad \;$ be a function $\; a \colon I \raro \Ad \;$ on $I$
with the values $\; a(i) \;,\; i \in I  \;,$ in the group $\Ad$ of
all permutations of $\; Y_d = \{1,2, \ldots ,d \} \;$. Consider a
$d$-extension $\; \bI_a \;$ generated by the function $a$ (See
Section \ref{ss4.2}).
We assume that the group $ \gG(a) $, generated by $\; a(i) , i \in
I, \;$ acts transitively on $Y_d$. This provides that the shift
$T_{\bI_a} $ and the skew product $\bT_{I,a}$ are ergodic.

We want to clarify: when is $\; \bI_a \;$ the canonical graph for
the corresponding Markov shift $T_{\bI_a}$ (Definition \ref{def
canon form}). Let $\pi \colon \bI_a \raro I$ be the projection and
$ (\pi,\psi) \in \Ext $. Since every homomorphism $ \psi \colon I
\raro I $ is an automorphism, the pair $(\pi,\psi)$ is
irreducible. Therefore $\; \bI_a \;$ is a the canonical graph iff
$\; d(T_{\bI_a}) = d \;$.
\begin{proposition}\label{d = d(T bI a)}
If the function $a$ satisfies the following condition
 \begin{equation}\label{modular a}
  \rho(i) = \rho(i') \;\; \RRaro \;\; a(i) = a(i')
  \;\;,\;\; i,i' \in I
 \end{equation}
then $\; d(T_{\bI_a}) = d \;$.
\end{proposition}
\proof Suppose the condition (\ref{modular a}) holds.
The Markov shift $T_{\bI_a} $ is isomorphic to the skew product
$\bT = \bT_{\rho,a}$, which acts on $X_\rho \times Y_d$ by
(\ref{bT rho a}). So that we have $d(T_{\bI_a}) = d(\bT) $ and by
Theorem \ref{simp mar} $\; d(\bT) \;=\; d_{\gga \colon \gb}(\bT)
\;$.

A direct computation, using (\ref{modular a}), the definition of
$\gga(T)$ and $\gb(\bT)$ and Proposition \ref{ga(T rho)}, shows that
$$
 \gb(\bT) = \gga(T_\rho) \times \gep_{Y_d}  \;\;\;,\;\;\;
 \gga(\bT) = \gga(T_\rho) \times \nu_{Y_d} \;.
$$
This means that any element of $\gga(\bT)$ consists precisely of
$d$ elements of the partition $\gb(\bT)$. By the definition of the
index $ d_{\gga \colon \gb}(\bT) $ we have $\; d_{\gga \colon
\gb}(\bT) = d \;$. Thus $\; d(T_{\bI_a}) = d \;$. \qed

\medskip

Taking into account Theorem \ref{th rho cohom} we have

\begin{corollary}\label{cor modular a}
Let $\; \pi_k \colon \bI_{a_k} \raro  I \;,\; k = 1,2, \;$ be two
$d$-extensions of the Bernoulli graph $(I,\rho)$, generated by
functions $ a_k \colon I \raro \Ad $, respectively, and suppose
both the functions $\; a_k \;,\;k=1,2 \;$ satisfy the condition
\ref{modular a}.
Then the Markov shifts $T_{\bI_{a_1}}$ and $T_{\bI_{a_2}}$ are
isomorphic iff $ a_1 $ and $ a_2 $ are conjugate in $\Ad$, i.e.
there exists $w_0 \in \Ad$ such that $\; a_2(i) \cdot w_0 = w_0
\cdot a_1(i) \;,\; i \in I \;$.
\end{corollary}
\begin{remark}\label{rem modular a}
It can be proved that for $d$-extension $\bI_a$, the condition
\ref{modular a} is equivalent to $\; d(T_{\bI_a}) = d \;$.
\end{remark}

\bigskip

\noindent {\bf Absolutely non-homogeneous $\bldrho$.}  Consider
the case , when $\rho$ is absolutely non-homogeneous (see Section
\ref{ss2.4}). this means that $\; \rho(i) \neq \rho(i') \;$ for
all $\;i \neq i' \;$ from $I$, i.e. the Bernoulli graph $(I,\rho)$
has no congruent edges.

In this case for any $\rho$-uniform graph $G$ there exists a {\bf
unique} homomorphism $\; \phi \colon G \raro I \;$. Therefore
Theorem \ref{phi bar d(T)} can be sharpened as follows

\begin{theorem}\label{phi bar d(T) sharp}
Let $G$ be a $\rho$-uniform  stochastic graph, which is
irreducible and satisfies the positive recurrence condition.
Suppose that $\rho$ is absolutely non-homogeneous. Then there
exist a unique homomorphism $\; \phi \in \Hom (G,I) \;$ and a
commutative diagram
\begin{equation}\label{diag bar H G}
 \xymatrix{   \bH \ar[d]^{\pi} \ar[r]^{\bpsi} & G \ar[d]^{\phi} \\
                               H \ar[r]^{\psi} & I   }
\end{equation}
such that
 \begin{enumerate}
  \item[(i)] The pair $\; (\pi,\psi) \in \Ext$ is a
  $(\pi,\psi)$-extension.
  \item[(ii)] $\; d = d(\phi) = d(T_G) \;$,
 \end{enumerate}
\end{theorem}
A natural question, which is arisen in connection with the
previous theorem is:
\begin{question}[{\bf Generalized Road Coloring Problem}]
\label{gen road prob}
Does Theorem \ref{phi bar d(T)} hold with $n=1$ in general case,
when $\rho$ is not necessarily absolutely non-homogeneous, i.e.
when $(I,\rho)$ has congruent edges ?
\end{question}
As we know, the problem is open even in the case, when the graph
$G$ is finite (See \cite{AsMaTu} and references therein.)

\bigskip

\noindent {\bf Homogeneous $\bldrho$ and Road Problem } Consider a
special case, when $\rho$ is homogeneous, i.e. $\; \rho(i) =
l^{-1} \;,\; i \in I \;$ with an integer $l = |I| \in \N$. Theorem
\ref{simp mar} and arguments adduced in Section \ref{ss2.4} imply
\begin{theorem}\label{} Suppose $\rho$ is homogeneous.
Then every ergodic $\rho$-uniform Markov shift $\; T_G \;$ is
isomorphic to a direct product $\; T_\rho \times \gs_d \;$ of the
Bernoulli shift $T_\rho$ and a cyclic permutation $\gs_d$ of
$Y_d$, where $d$ is the period of $ T_G $. If, in addition, $ T_G
$ is exact, then it is isomorphic to the Bernoulli shift $T_\rho$,
herewith, there exists $n \in \N$ and a homomorphism $\; \phi
\colon G\pn \raro I \;$ of degree $1$.
\end{theorem}
The result was proved earlier in \cite{Ru$_3$} for finite $G$ and
in \cite{Ru$_6$} for general case.

If $G$ is finite and $\rho$ is homogeneous Question \ref{gen road
prob} is a reformulation of well-known Road Coloring Problem (See
\cite{Fr}, \cite{O'B}, \cite{AdGoWe}, \cite{Ki}). As we know, the
problem is still open.


\bigskip

\subsection{Some (p,q)-uniform graphs.}
\label{ss5.4} We construct some simple examples to illustrate the
case, when the $\psi$-part in the canonical pair $(\pi,\psi)$ is
not trivial.

Let $\; I = \{0,1\}\;$ and  $\; \rho = (p,q) \;$, where $\; 0 < p < 1 \;$
and $\; q=1-p \;$.
Given $\; n \in \N \;$ consider the following random walk on
$\; J_n := \{ 1,2, \ldots ,n \} \;$
\begin{equation}\label{FDR}
   \xymatrix@C=3pc{
            1       \ar@(ul,dl) []_{q}   \ar@/_/ [r]_{p}
    &       2       \ar@/_/ [l]_{q}      \ar@/_/ [r]_{p}
    &       \;\;\;\;\hdots\;\;\;\;
                    \ar@/_/ [l]_{q}      \ar@/_/ [r]_{p}
    &       n       \ar@/_/ [l]_{q}      \ar@(dr,ur) []_{p}    } \;,
\end{equation}
which is known as a {\bf Finite Drunkard Ruin}.
We set here: $\; H:=I \times J_n \;$ , $\; H\po0 := J_n \;$ and
$$
  s(h) = j \;,\; t(h) = f_ij \;,\; \psi (h) = i
  \;\;,\;\; h = (i,j) \in H \;,
$$
where the maps $\; f_i \colon J_n/ \raro J_n \;,\; i = 0,1 \;,$
are defined by
$$
  f_1 j = \min{(j+1,n)} \;\;,\;\;  f_0 j = \max{(j-1,1)}
   \;\;\;,\;\;\;  j \in J_n
$$
and the weights of edges $\; p(h) \;,\; h \in H \;$ are given
according to (\ref{FDR}) by $\; p(1,j) = p \;$, $\; p(0,j) = q
\;$.

Then the finite stochastic graph $H$ is irreducible and
$\rho$-uniform, $\; \psi \in \Hom (H,I) \;$. The semigroup $\cS
(\psi)$, generated by $\{f_0,f_1\}$, is $1$-contractive, since $\;
(f_0)^n (J_n) = \{1\} \;$. Whence, $\; d(\psi) = 1 \;$ and the
Markov shift $T_H$ is isomorphic to the Bernoulli shift $T_\psi$.

\medskip

Given $p$ and $n$ we construct a $\Z_2$-extension $\bH_a$ of the
graph $H$, where $a : H \raro \Z_2$ and $\Z_2 := \{0,1\}$ be the
cyclic group of order $2$.

Define $\; a : H=I \times J_n \ni h=(i,j) \raro a(h) \in \Z_2 \;$ by
\begin{equation}\label{a(h)}
a(i,j) \;=\;            \left\{
     \begin{array}{ll}
         1 \;\;,\;\; & if \;\; (i,j)    \;=\;   (1,1) \; \\
         0 \;\;,\;\; & if \;\; (i,j)   \;\neq\; (1,1) \;.
     \end{array}
                        \right.
\end{equation}
Then the corresponding graph $\bH_a$ has the form
\begin{equation}\label{FDR2}
   \xymatrix@C=3pc{
      **[r] 11      \ar@(ul,dl) []_{q}  \ar [rd]_(.75){p}
    &       21      \ar@/_/ [l]_{q}     \ar@/_/ [r]_{p}
    &       \;\;\;\hdots\;\;\;
                    \ar@/_/ [l]_{q}     \ar@/_/ [r]_{p}
    & **[l] n1      \ar@/_/ [l]_{q}     \ar@(dr,ur) []_{p}
    &      z=1                                             \\
      **[r] 10      \ar@(dl,ul) []^{q}  \ar [ru]^(.75){p}
    &       20      \ar@/^/ [l]^{q}     \ar@/^/ [r]^{p}
    &       \;\;\;\hdots\;\;\;
                    \ar@/^/ [l]^{q}     \ar@/^/ [r]^{p}
    & **[l] n0      \ar@/^/ [l]^{q}     \ar@(ur,dr) []^{p}
    &      z=0                                             }
\end{equation}
for $\; n > 2 \;.$ and
\begin{equation}\label{FDR3}
   \xymatrix@C=3pc{
        **[r] 11      \ar@(ul,dl) []_{q}  \ar@/^3pc/ [d]^{p}
    && **[r] 11      \ar@(ul,dl) []_{q}  \ar [rd]_(.75){p}
    &   **[l] 21      \ar@/_/ [l]_{q}     \ar@(dr,ur) []_{p}
    &        z=1                                              \\
        **[r] 10      \ar@(ul,dl) []_{q}  \ar@/_2pc/ [u]^{p}
    && **[r] 10      \ar@(dl,ul) []^{q}  \ar [ru]^(.75){p}
    &   **[l] 20      \ar@/^/ [l]^{q}     \ar@(ur,dr) []^{p}
    &        z=0                                               }
\end{equation}
for two special cases $\; n = 1,2 \;$

\medskip

Suppose $\; p \neq q \;$. We claim in this case that for all $n
\in \N $ the graphs (\ref{FDR3}) and (\ref{FDR2}) are canonical.
Indeed, $\; d(\pi_H) = d(T_{\bH_a}) = 2 \;$, since $ \rho = (p,q)
$ is absolutely non-homogeneous. In order to check the
irreducibility of the $2$-extension $\; (\pi_H,\psi) \colon \bH_a
\raro H \raro I \;$ consider the semigroup $\bcS = \bcS(\pi,\psi)$
and its persistent partitions $\cR (\bcS)$.

The semigroup $\bcS$ is generated by $\; \{\baf_i , i\in I \}$, where
$$
 \baf_i \; (j\;,\;z) \;=\; (f_i \; j \;,\; z + a(i,j) \pmod 2 ) \;\;,\;\;
       (j,z) \in J_n \times \Z_2 \;.
$$
A direct computation shows that for $\; n = 1,2 \;$ any transversal
partition of $J_n \times \Z_2 $ is persistent in the sense of
Definition \ref{pers part} , and for $\; n > 2 \;$ there exists a
non-persistent transversal partition. Naimly, the partition,
consisting of two sets of the following "alternating" form
$$
           \{  (1,0) , (2,1) , (3,0) , (4,1) , \ldots  \}
 \;\;,\;\; \{  (1,1) , (2,0) , (3,1) , (4,0) , \ldots  \} \;,
$$
is so.
Moreover, this is the only transversal partition, which is
not persistent. This implies that for every $ n \in \N$ the
persistent partitions $\cR (\bcS)$ separate points of $ J_n \times
\Z_2 $ in the sense of Remark \ref{irred ext} and the
$2$-extension $(\pi_H, \psi)$ is irreducible. Thus
\begin{itemize}
 \item For all $n \in \N $ and $p \neq q$ the graphs $\bH_a$ are canonical
 graphs for the corresponding shifts $T_{\bH_a}$.
\end{itemize}
Just in the same way we can consider the following
{\bf Infinite Drunkard Ruin}
\begin{equation}\label{FDR4}
   \xymatrix@C=3pc{
      **[r] 1      \ar@(ul,dl) []_{q}  \ar@/_/ [r]_{p}
    &       2      \ar@/_/ [l]_{q}     \ar@/_/ [r]_{p}
    &                   \;\;\;\hdots\;\;\;
                   \ar@/_/ [l]_{q}     \ar@/_/ [r]_{p}
    &       n      \ar@/_/ [l]_{q}     \ar@/_/ [r]_{p}
    &              \ar@/_/ [l]_{q}     \;\;\;\hdots      }
\end{equation}
where  $\; H:=I \times \N \;$ , $\; H\po0 := \N \;$.

Suppose $\; p < q \;$. Then the corresponding Markov chain is
positively recurrent and the Markov shift $T_H$ is isomorphic to
the Bernoulli shift $T\rho$.

Again define the functions $\; a : H=I \times \N \ni h=(i,j) \raro
a(h) \in \Z_2 \;$ by (\ref{a(h)}). Then $\Z_2$-extension $\bH_a$
of the graph $H$ (\ref{FDR5}) has the form
\begin{equation}\label{FDR5}
   \xymatrix@C=3pc{
      **[r] 11               \ar@(ul,dl) []_{q}  \ar [rd]_(.75){p}
  &         21               \ar@/_/ [l]_{q}     \ar@/_/ [r]_{p}
  &   \;\;\;\hdots\;\;\; \ar@/_/ [l]_{q}     \ar@/_/ [r]_{p}
  &         n1               \ar@/_/ [l]_{q}     \ar@/_/ [r]_{p}
  &                          \ar@/_/ [l]_{q}     \;\;\;\hdots      \\
      **[r] 10               \ar@(dl,ul) []^{q}  \ar [ru]^(.75){p}
  &       20                 \ar@/^/ [l]^{q}     \ar@/^/ [r]^{p}
  &   \;\;\;\hdots\;\;\; \ar@/^/ [l]^{q}     \ar@/^/ [r]^{p}
  &       n0                 \ar@/^/ [l]^{q}     \ar@/^/ [r]^{p}
  &                          \ar@/^/ [l]^{q}     \;\;\;\hdots       }
\end{equation}
It can be shown in this case that any transversal set is
persistent. Thus
\begin{itemize}
 \item If $p < q$ the graph $\bH_a$ (\ref{FDR5}) is the canonical graph
  for the shift $T_{\bH_a}$.
\end{itemize}
Note that the shift $T_{\bH_a}$ is a $\Z_2$-extension of the
Bernoulli shift $T_{p,q}$, therefore, $T_{\bH_a}$ has a
$4$-element one-sided generator. On the other hand the shift is
not isomorphic to Markov shifts on finite state spaces. Thus
\begin{itemize}
\item If $p < q$ the one-sided Markov shift $T_{\bH_a}$ has no finite
one-sided Markov generator.
\end{itemize}


\end{document}